\documentclass[reqno,draft]{amsart}

\theoremstyle{plain}
\newtheorem{lemma}{Lemma}
\newtheorem{theorem}[lemma]{Theorem}
\newtheorem{corollary}[lemma]{Corollary}

\newtheorem{proposition}[lemma]{Proposition}
\newtheorem{definition}[lemma]{Definition}

\theoremstyle{remark}

\newcommand*  {\R} {{\mathbb R}}

\newcommand*{\ang}[1]{\left\langle #1 \right\rangle}

\def\ee{\epsilon}
\def\aa{\alpha}

\def\dd{\delta}

\def\ss{\sigma}
\def\tt{\theta}

\def\ll{\lambda}
\def\la{\lambda}
\def\Dd{\Delta}
\def\Gg{\Gamma}
\def\Om{\Omega}

\def\pp{\partial}

\begin{document}

\title[PG Equations with Hyper--Diffusion]
{A ``horizontal" hyper--diffusion $3-D$ thermocline 
planetary geostrophic model: well-posedness and long time behavior}

\date{June 2, 2003}

\author[C. Cao]{Chongsheng Cao} 
\address[C. Cao]
{Department of Mathematics  \\
University of Nebraska-Lincoln \\
Lincoln, NE 68588-0323  \\
USA.}
\email{ccao@math.unl.edu}

\author[E.S. Titi]{Edriss S. Titi}
\address[E.S. Titi]
{Department of Mathematics \\
and  Department of Mechanical and  Aerospace Engineering \\
University of California \\
Irvine, CA  92697-3875, USA \\
{\bf Also}  \\
Department of Computer Science and Applied Mathematics \\
Weizmann Institute of Science  \\
Rehovot 76100, Israel.}  
\email{etiti@math.uci.edu}

\author[M. Ziane]{Mohammed Ziane}
\address[M. Ziane]
{Department of Mathematics \\
University of Southern California \\
Los Angeles, CA 90089-1113, 
USA.}
\email{ziane@math.usc.edu}

\begin{abstract}
In this paper we study a three dimensional thermocline planetary geostrophic 
``horizontal" hyper--diffusion model of the gyre-scale midlatitude ocean. 
We show the global existence and uniqueness of the weak and strong 
solutions to this model.
Moreover, we establish the existence of a finite dimensional global
attractor to this dissipative evolution system. Preliminary computational
tests indicate that our hyper--diffusion model does not exhibit any of 
the nonphysical instabilities near the literal boundary which are 
observed numerically in other models.
\end{abstract}

\maketitle


\section{Introduction}   \label{S-1}
The planetary geostrophic (PG) models, the adiabatic and inviscid
form of ``thermocline'' or ``Phillips type II'' equations,
 are derived by standard
scaling analysis for gyre-scale oceanic motion (see \cite{PJ84}, \cite{Ph}, 
\cite{RS59} and \cite{WE59}).  In their simplest
dimensionless $\beta-$plan form, these equations are:

\begin{eqnarray}
&&\hskip-.8in
\partial_x p - f  v_2  = 0  \label{PG-1}  \\
&&\hskip-.8in
\partial_y p + f  v_1 = 0  \label{PG-2}  \\
&&\hskip-.8in
\partial_z p  +T =0    \label{PG-3}  \\
&&\hskip-.8in
\nabla\cdot v+\partial_z  w =0   \label{PG-4} \\
&&\hskip-.8in
\partial_t T + v\cdot \nabla T + w \partial_zT = 0 \, ,  \label{PG-5} 
\end{eqnarray}
in the domain $\Omega = \{ (x,y,z) : 
(x,y) \in M \subset \R^2,~ {\mbox{and}}~z \in (-h, 0) \}$. 
Here $\nabla =(\pp_x, \pp_y)$, 
$v=(v_1,v_2)$ denotes the horizontal velocity field, 
$w$ the vertical velocity,
 $p$ is the pressure,
 $T$ is the temperature, and $ f = f_0 + \beta y$ is the Coriolis
 parameter.  A diffusive term, $\kappa_v \partial_z^2T$, is commonly added 
to equation~(\ref{PG-5}) as a leading order approximation to the
effect of microscale turbulent mixing. Thus equation~(\ref{PG-5})
becomes:
\begin{equation}
\partial_t T +  v\cdot\nabla T + w \partial_zT  =  \kappa_v \partial_z^2T .  
\label{PG-6} 
\end{equation}

In \cite{SV97A} it is  argued, based on physical
grounds, that in closed ocean basin and with the no-normal-flow
boundary conditions,  this model can be solved  only
in restricted domains which are bounded away from the lateral
boundaries, $\partial M \times (-h,0)$. Thus, it cannot 
be utilized in the study of the
large-scale circulation. Furthermore, it has
been  pointed out, numerically, by
   \cite{CV86} that arbitrarily  short     linear 
disturbances (disturbances that are supported at small spatial scales)
will grow arbitrarily fast when the flow becomes
baroclinically unstable. This nonphysical growth  at small scales
is a signature of mathematical    ill-posedness of this model near
unstable baroclinic mode. Therefore, Samelson and Vallis proposed
in \cite{SV97A} various  simple friction and 
diffusion schemes to overcome these physical and numerical
 difficulties. In particular, they propose a linear drag
(or viscosity)  in the horizontal momentum equations and a
horizontal diffusion in the thermodynamic 
equation~(\ref{PG-6}). Thus the full dimensionless system becomes:

\begin{eqnarray}
&&\hskip-.8in
\partial_x p - f  v_2  = -\epsilon v_1+ \nu_v\partial_z^2 v_1
+\nu_h  \Delta_h v_1 
\label{PG-11a}  \\
&&\hskip-.8in
\partial_y p + f  v_1  = -\epsilon v_2 +\nu_v \partial_z^2v_2+\nu_h  
\Delta_h v_2\label{PG-12a}  \\
&&\hskip-.8in
\partial_z p   + T =0    \label{PG-13a}  \\
&&\hskip-.8in
\nabla v+\partial_z w =0   \label{PG-14a} \\
&&\hskip-.8in
\pp_t T + v\cdot \nabla T + w\partial_z T = 
 \kappa_v \partial_z^2T+ \kappa_h \Delta_h T \, , \label{PG-15a}
\end{eqnarray}
where $\Delta_h=\frac{\pp^2}{\pp x^2}+ \frac{\pp^2}{\pp y^2} $. 
Notice that the incompressibility and hydrostatic balance are
retained in the above system. Here, the coefficients $\nu_h$, 
 $\nu_v$, $ \epsilon$ and $ \kappa_h$ are non-negative and 
small. 

\medskip

In the case where $\nu_h > 0$, 
 $\nu_v > 0$, $ \epsilon \ge 0 $,  $ \kappa_v > 0$  and $ \kappa_h > 0$, 
with the ocean being driven by the wind stresses at the top
surface, and with  no-slip boundary conditions
 and  no-heat fluxes on the side walls, $\Gamma_s=\partial M \times (-h,0)$,
 and at  the bottom, $M \times \{z= -h\}$,
 the above system has been studied analytically by
 \cite{STW98} and  \cite{STW00}.
Under this type of boundary conditions the first two authors
have improved in  \cite{CT01}  the results of 
 \cite{STW98} and  \cite{STW00}
and established global regularity and well-posedness to the
 system~(\ref{PG-11a})-(\ref{PG-15a}).
 In particular, they show in \cite{CT01}  the existence
of finite dimensional global attractor to this viscous  three dimensional
planetary geostrophic (PG)  model.  This global well-posedness result  provides
a rigorous justification to the scheme suggested by Samelson
and Vallis  \cite{SV97A} under the above conditions
on the coefficients. 

\medskip

On the other hand, the case where  $\nu_h = 0$, 
 $\nu_v = 0$, $ \epsilon > 0 $,  $ \kappa_v > 0$ and $ \kappa_h > 0$,
has been suggested in related schemes by other authors
 \cite{Killworth},
\cite{Salmon-86}, \cite{Salmon-90}, 
 \cite{Winton-Sarachik}, and \cite{Zhang}. Equations ~(\ref{PG-11a})
-~(\ref{PG-15a}) become
\begin{eqnarray}
&&\hskip-.8in
\partial_x p - f  v_2  = -\epsilon v_1
 \label{PG-11}  \\
&&\hskip-.8in
\partial_y p + f  v_1  = -\epsilon v_2  \label{PG-12} \\
&&\hskip-.8in
\partial_z p   + T =0    \label{PG-13}  \\
&&\hskip-.8in
\nabla v+\partial_z w =0   \label{PG-14} \\
&&\hskip-.8in
\pp_t T + v\cdot \nabla T + w\partial_z T = 
 \kappa_v \partial_z^2T+ \kappa_h \Delta_h T \, .\label{PG-15}
\end{eqnarray}
 In this situation, and by 
following   \cite{SV97A},
the frictional-geostrophic relations (\ref{PG-11}) and (\ref{PG-12}) 
can be solved locally for the horizontal velocities to give:
\begin{equation}
\label{HV}
v_1=-\gamma (\epsilon p_x +f p_y),
 \qquad v_2=\gamma (f p_x - \epsilon p_y) \, ,
\end{equation}
where $ \gamma= (f^2 +\epsilon^2)^{-1}$. Assuming we have
a nice solution up to the boundary, equations~(\ref{HV}) 
and~(\ref{PG-13}) imply:
\begin{equation}
\label{HVZ}
\partial_z v_1=-\gamma (\epsilon T_x +f T_y),
 \qquad \partial _z v_2=\gamma (f T_x - \epsilon T_y) \, .
\end{equation}
As a result of the above, the no-normal-flow boundary condition
on the lateral boundary yields  the following boundary condition
on the temperature:
\begin{equation}
\label{BCT1}
\epsilon \frac {\pp T}{\pp n} +f \frac{\pp T}{\pp s} = 0 \qquad 
{\mbox{on the lateral boundary}},~ \Gamma_s=\partial M \times (-h,0)\, ,
\end{equation}
where  $ \frac{\pp T}{\pp n}$ and $ \frac{\pp T}{\pp s}$ denote the normal
and right-hand tangential derivative, respectively, so that with 
$\vec{e} = \left( \ee n_1 -f n_2, f n_1 +\ee  n_2
\right)/(\ee^2+f^2)^{1/2}$, we have the following oblique boundary condition
on the temperature.
\begin{eqnarray}
&&\hskip-.8in
\frac{\pp T}{\pp \vec{e}} =0,  \quad \mbox{on } \Gg_s.
\label{XBBB0}
\end{eqnarray}
On the other
hand,  the no heat flux at the lateral boundary gives:
\begin{equation}
\label{BCT2}
-\kappa_h \frac {\pp T}{\pp n} =0 \qquad 
{\mbox{on the lateral boundary}},~ \partial M \times (-h,0)\, .
\end{equation}
Based on the above it is argued in  \cite{SV97A} that
in the presence of rotation, i.e., $f\not= 0$, and in 
order to be able to  satisfy both boundary conditions~(\ref{BCT1}) 
and~(\ref{BCT2}) one has to add to  the thermodynamics
 equation a higher
order (biharmonic) horizontal  diffusion. Otherwise, the 
 problem~(\ref{PG-11})-(\ref{PG-15}) subject to the additional
 boundary
conditions~(\ref{BCT1}) 
and~(\ref{BCT2})
is over determined and hence is ill-posed, which explains 
the cause for the observed numerical instabilities near the 
lateral boundary. In order to resolve this  discrepancy
 Samelson and Vallis propose in \cite{SV97A} 
the following ``horizontal hyper--diffusion'' planetary geostrophic (PG) model:
\begin{eqnarray}
&&\hskip-.8in
 \partial_x p- f  v_2  = -\epsilon v_1 \label{PG-21}  \\
&&\hskip-.8in
\partial_y p + f  v_1  = -\epsilon v_2 \label{PG-22}  \\
&&\hskip-.8in
\partial_z p   + T =0    \label{PG-23}  \\
&&\hskip-.8in
\nabla \cdot v +\partial_z w =0   \label{PG-24} \\
&&\hskip-.8in
\pp_t T + v\cdot \nabla T + w\partial_z T = 
 \kappa_v \partial_z^2T+ \kappa_h \Delta_h T - \lambda\Delta_h^2 
 T    \, , \label{PG-25}
\end{eqnarray}
 subject to the no-normal-flow together with the
 boundary conditions given in ~(\ref{BCT1}) on the lateral boundary. 
However, the no-heat flux boundary condition~(\ref{BCT2}) is replaced 
by the boundary condition:
\begin{equation}
\label{BCT3}
-\kappa_h \frac {\pp T}{\pp n} +
 \lambda \Delta_h \frac {\pp T}{\pp n}   =0 \qquad 
{\mbox{on the lateral boundary}},~\Gamma_s= \partial M \times (-h,0).
\end{equation}
It is worth stressing that the differences between the boundary conditions 
(\ref{BCT1}) and (\ref{BCT2}) is due to the Coriolis parameter. 
Therefore, it is natural to require $\la \rightarrow 0$ when 
$f \rightarrow 0$. Since the diffusion operator 
$\kappa_v \partial_z^2T+ \kappa_h \Delta_h T - \lambda\Delta_h^2 T$ in the model
(\ref{PG-21})--(\ref{PG-25}) with the boundary conditions (\ref{XBBB0}) and
(\ref{BCT3}) may not positive, the numerical instabilities have been observed 
near the lateral boundary.

Motivated by \cite{SV97A} we propose below a variant of the system 
(\ref{EQ-1})--(\ref{BT-3}) which is globally well-posed and which
possesses a finite dimensional global attractor.

\bigskip

Let 
$\Om= M \times (-h, 0)\in \mathbb{R}^3,$
where $M$ is a smooth domain in $\mathbb{R}^2,$  or 
$M=(0,1) \times (0,1),$ the PG equations with
friction and diffusion can be written as:
\begin{eqnarray}
&&\hskip-.8in
\nabla p + f \vec{k} \times v +\mathcal{D} = 0  \label{G-EQ-1}  \\
&&\hskip-.8in
\pp_z p  + T =0    \label{G-EQ-2}  \\
&&\hskip-.8in
\nabla \cdot v +\pp_z w =0   \label{G-EQ-3} \\
&&\hskip-.8in
\pp_t T  + v \cdot \nabla T + w \pp_z T = 
K_v T_{zz} - \nabla \cdot q(T),  \label{G-EQ-4} 
\end{eqnarray}
where $\mathcal{D}$ is the friction or dissipation of momentum 
and $\nabla \cdot q(T)-K_v T_{zz}$ is the heat diffusion.

Naturally, there are two friction schemes: one is the linear drag, i.e.,   
$\mathcal{D} = \ee v;$
the other is the conventional eddy viscosity, i.e.,
$
\mathcal{D} = -\ee (A_h \Dd v + A_v \pp_{zz}v),
$
where $A_h$ and $A_v$ are small positive constants. 
The PG model with conventional eddy viscosity 
and the diffusion
$q(T)=-K_h \nabla T,$
has been studied in {\bf \cite{CT01}, \cite{STW98}, \cite{STW00}}. 
Here we concentrate on the case of linear drag.
As explained above,  in a closed basin and with a linear drag scheme, 
it is appropriate that $\nabla \cdot q$ be a fourth--order diffusion. 
Denoting, from now on, by $\Dd = \Dd_h=\pp_x^2 +\pp_y^2.$ 
Our goal in this paper is to study the model with $\mathcal D= \ee v$, 
and with ``horizontal" hyper--diffusion 
\begin{eqnarray}
&&\hskip-.18in 
q(T) =\la H \nabla ( \nabla \cdot (H^T \nabla T)) 
- K_h \nabla T + 
\mu \nabla T_{zz},   \label{QQ} 
\end{eqnarray} 
where 
\begin{eqnarray}
&&\hskip-.18in 
H = \left( \begin{array}{ll}
1 & -f/\ee \\ f/\ee & 1 \end{array} \right),   \label{HH} 
\end{eqnarray} 
and $H^T$ is the transpose of $H$, and $\la$, $\mu$ and $\ee$ are 
small positive constants.
That is we study the following model:
\begin{eqnarray}
&&\hskip-.8in
\nabla p + f \vec{k} \times v +\ee v = 0  \label{EQ-1}  \\
&&\hskip-.8in
\pp_z p  + T =0    \label{EQ-2}  \\
&&\hskip-.8in
\nabla \cdot v +\pp_z w =0   \label{EQ-3} \\
&&\hskip-.8in
\pp_t T  + v \cdot \nabla T + w \pp_z T =Q+ 
K_v T_{zz} 
-\nabla \cdot \left( \la H \nabla ( \nabla \cdot (H^T \nabla T)) 
- K_h \nabla T + 
\mu \nabla T_{zz} \right),  \label{EQ-4} 
\end{eqnarray}
where $Q$ is a given heat source.

\medskip

Now we impose the  appropriate boundary conditions to 
this friction and hyper--diffusion PG model (\ref{EQ-1})--(\ref{EQ-4}). 
The natural boundary conditions are  
 no--normal flow condition for the velocity field $(v,w)$ 
and non--flux boundary condition for the temperature $T$
(see, e.g., {\bf \cite{PJ87}, \cite{SV97}, \cite{SV97A}, \cite{SD96}}): 
\begin{eqnarray}
&&\hskip-.8in
\mbox{on } \Gg_u: w=0;    \label{B-1}\\
&&\hskip-.8in
\mbox{on } \Gg_b:  w=0;  \label{B-2}\\
&&\hskip-.8in
\mbox{on } \Gg_s: 
v \cdot \vec{n} =0, 
\label{B-3}  \\
&&\hskip-.8in
\mbox{on } \Gg_u:   -K_v \frac{\pp T }{\pp z} = \aa (T-T^*);    
\label{BT-1}\\
&&\hskip-.8in
\mbox{on } \Gg_b: \frac{\pp T }{\pp z} = 0; 
\label{BT-2}\\
&&\hskip-.8in
\mbox{on } \Gg_s: q (T) \cdot \vec{n}  =0;    \label{BT-3}
\end{eqnarray}
where $\Gg_u, \Gg_b$ and $\Gg_s$ denote the boundary of $\Om$ defined as:
\begin{eqnarray}
&&\hskip-.8in
\Gg_u = \{ (x,y,z)  \in \Om : z=0 \},  \\
&&\hskip-.8in
\Gg_b = \{ (x,y,z)  \in \Om : z=-h \},  \\
&&\hskip-.8in
\Gg_s = \{ (x,y,z)  \in \Om : (x,y) \in \partial M  \},  
\end{eqnarray}
$\aa$ is a positive constant, $\vec{n}=(n_1,n_2)$ is the normal
vector of $\Gg_s$,  and $T^* (x,y)$ is a typical top surface 
temperature profile.  Furthermore,
the no--normal flow condition implies a boundary 
condition on the temperature, as explained in (\ref{HV}), (\ref{HVZ}), 
and (\ref{BCT1}):
\begin{eqnarray}
&&\hskip-.8in
 \frac{\pp T}{\pp \vec{e}} =0, \quad \mbox{on } \Gg_s,
\label{XB}
\end{eqnarray}
where $\vec{e} = \left( \ee n_1 -f n_2, f n_1 +\ee  n_2
\right)/(\ee^2+f^2)^{1/2}.  $
Therefore, no--normal flow (\ref{XB}) and no heat--flux boundary conditions
(\ref{BT-1})--(\ref{BT-3}) are natural and proper boundary conditions for 
the PG model (\ref{EQ-1})--(\ref{EQ-4}). Since (\ref{XB}) is in fact a boundary 
condition on the temperature, thus, any additional and incompatible boundary 
condition to the temperature, on top of  (\ref{BT-1})--(\ref{BT-3})
would make the fourth--order diffusion $\nabla \cdot q$ overdetermined. 
Finally, the model is supplemented with the initial condition:
\begin{eqnarray}
&&\hskip-.18in 
T(x,y,z,0) = T_0(x,y,z), \label{INIT} 
\end{eqnarray} 
where $T_0$ is a given function.

We observe that in the case when $\beta =0$, i.e.~$f=f_0$, our
``horizontal''  hyper-diffusion  term $\nabla \cdot q(T)$ 
reduces to the form
 $ \lambda\Delta^2 T - \mu \Delta \partial_z^2T- \kappa_h \Delta T $,
which is in the spirit of the hyper-diffusion term in~(\ref{PG-25})
that was proposed by Samelson and Vallis in \cite{SV97} and \cite{SV97A}.
However, due to the fact that $\beta \neq 0$ our proposed
hyper-diffusion
term takes a more involved form, which is necessary to 
guarantee the dissipativity of this operator 
 under the given physical boundary
 conditions~(\ref{B-1})-(\ref{BT-3}). It is worth adding, that the
 present PG formulation (\ref{EQ-1})--(\ref{BT-3})
has been explored with
some preliminary computations by Samelson \cite{SAMELSON} 
using a modified version of the
$\beta-$plane numerical code developed in  \cite{SV97A}.  
In the $\beta-$plane case (for which
the horizontal coordinates are Cartesian $x$ and $y$, and the
Coriolis parameter $f=f_0+\beta y$, with $f_0$ and $\beta$
constant), the additional horizontal diffusion terms in the
thermodynamic equation reduce to a single term proportional
to $T_{xx}$.  The preliminary computations indicate
that inclusion of these additional terms only slightly
modifies the previous numerical solutions reported in \cite{SV97A}
which uses hyper-diffusion term of the form suggested
in~(\ref{PG-25}).
With the rigorous analytical  results
proved here for the modified system(\ref{EQ-1})--(\ref{BT-3}),
 this provides new
theoretical and mathematical support for the approach and
analysis of \cite{SV97}.  Note that these
preliminary computations did not include the additional
mixed horizontal--vertical diffusion term $\mu \Delta T_{zz}$,
which is a crucial term for our rigorous mathematical analysis.

\section{Preliminaries}

It is natural to assume that $T^*$ satisfies the compatibility 
boundary  conditions:
\begin{eqnarray}
&&\hskip-1in 
\frac{\pp T^*}{\pp \vec{e}} = 0,   
\qquad \quad \quad \mbox{on } \pp M.     \label{TSTAR-1}  \\
&&\hskip-1in 
q(T^*) \cdot \vec{n}=0,   
\qquad \;  \mbox{on } \pp M.     \label{TSTAR-2}
\end{eqnarray}
Let  $ \widetilde{T}=T - T^*.$
Due to the compatibility 
boundary  conditions (\ref{TSTAR-1}) and (\ref{TSTAR-2}), 
it is clear that $\widetilde{T}$ satisfies the following 
homogeneous boundary conditions:
\begin{equation}
\left. \frac{\pp \widetilde{T}}{\pp z }
\right|_{z=-h}= 0; \quad 
\left. {\left( \frac{\pp \widetilde{T}}{\pp z} 
+ \frac{\aa}{K_v} \widetilde{T} \right) }
\right|_{z= 0}= 0;\quad 
\; \left. \frac{\pp \widetilde{T}}{\pp \vec{e}} \right|_{\Gg_s}= 0; 
\quad \; \left. q(\widetilde{T}) \cdot \vec{n} \right|_{\Gg_s}= 0. 
\label{HB}
\end{equation}
By replacing $T$ by $\widetilde{T}$, we have the following  new 
formulation of the system (\ref{EQ-1})--(\ref{INIT}):
\begin{eqnarray}
&&\hskip-.8in
\pp_t \widetilde{T}  + \nabla \cdot q(\widetilde{T}) 
- K_v \widetilde{T}_{zz} 
+ v \cdot \nabla \widetilde{T} 
+w \pp_z \widetilde{T}+   v  \cdot \nabla T^* =Q^*, 
 \label{TEQ-1}  \\ 
&&\hskip-.8in
\left. {\left( \pp_z \widetilde{T}+ \frac{\aa}{K_v}  \widetilde{T} 
 \right)} \right|_{z=0}= 0; \quad \left. \pp_z \widetilde{T}  
\right|_{z=-h}= 0; \quad \left. \frac{\pp \widetilde{T}}{\pp \vec{e}} 
\right|_{\Gg_s}= 0;
\quad \left. q(  \widetilde{T}) \cdot \vec{n} \right|_{\Gg_s}= 0,      
\label{TEQ-2}    \\
&&\hskip-.8in
\widetilde{T} (x,y,z,0) = T_0 (x,y,z)-T^*(x,y). \label{TEQ-3} 
\end{eqnarray}
{F}rom now on $q(\widetilde{T})$ is given by (\ref{QQ}).
Here $v$ and $w$ are determined by the use of (\ref{HVZ}) and 
(\ref{EQ-3}), and 
the fact that the average of $v$ in the $z$-direction is zero. 
This can be easily seen by 
integration of (\ref{EQ-1}) and (\ref{EQ-2}), which yields
 the following system:
\begin{eqnarray}
&&\hskip-.8in
\nabla \overline{p} + f \vec{k} \times \overline{v} + \ee \overline{v} = 0.
 \label{AV-EQ}  \\
&&\hskip-.8in
\nabla \cdot \overline{v} = 0,    \label{ADIV-1}   \\
&&\hskip-.8in
\overline{v} \cdot \vec{n} = 0.   \label{AV-BD-1}
\end{eqnarray}
Here $\overline{v}$ and $ \overline{p}$ are the averages of $v$ and 
$p$ in the 
$z$-direction.
Hence, multiplying by $\overline{v}$ and integrating over $M$, we 
obtain that 
$\overline{v}=0.$
The exact expressions of $v_1$, $v_2$ and $w$ in terms of $T$ are:
\begin{eqnarray}
&&\hskip-.8in
v_1 =  \int_{-h}^z 
\frac{\ee \widetilde{T}_x(x,y,\xi,t) 
+f \widetilde{T}_y(x,y,\xi,t)}{\ee^2 +f^2} d\xi +
\left(z+\frac{h}{2}\right) 
\frac{\ee T^*_x(x,y) +f T^*_y(x,y)}{\ee^2 +f^2}-  
\label{V1_N}        \\
&&\hskip-0.48in
- \frac{1}{h} \int_{-h}^0 \int_{-h}^{\eta} 
\frac{\ee \widetilde{T}_x(x,y,\xi,t) 
+f \widetilde{T}_y(x,y,\xi,t)}{\ee^2 +f^2} d\xi d \eta,
\nonumber        \\
&&\hskip-0.8in
v_2 =  \int_{-h}^z 
\frac{-f \widetilde{T}_x(x,y,\xi,t) 
+\ee \widetilde{T}_y(x,y,\xi,t)}{\ee^2 +f^2} d\xi 
+\left(z+\frac{h}{2}\right) 
\frac{-f T^*_x(x,y) +\ee T^*_y(x,y)}{\ee^2 +f^2} - 
\label{V2_N}         \\
&&\hskip-0.48in
- \frac{1}{h} \int_{-h}^0 \int_{-h}^{\eta} 
\frac{-f  \widetilde{T}_x(x,y,\xi,t) 
+\ee \widetilde{T}_y(x,y,\xi,t)}{\ee^2 +f^2} d\xi d\eta,
\nonumber   
\end{eqnarray}
Using (\ref{EQ-3}), (\ref{B-1}) and (\ref{B-2}) we have 
\begin{equation}
w(x,y,z)=-\int_{-h}^z (\nabla \cdot v) (x,y,\xi) d\xi.
\label{WWWW}
\end{equation}
From (\ref{V1_N}), (\ref{V2_N}) and (\ref{WWWW}) we obtain
\begin{eqnarray}
&&\hskip-0.8in
w =  - \int_{-h}^z \int_{-h}^{\eta} \left[
\frac{\ee \Dd \widetilde{T} (x,y,\xi,t) 
- f_0 \widetilde{T}_x(x,y,\xi,t)}{\ee^2 +f^2} +
\right. 
\label{W_N}          \\
&&\hskip-0.48in 
\left. + \frac{2 f_0 f}{(\ee^2 +f^2)^2} 
{\left( -f \widetilde{T}_x(x,y,\xi,t) 
+\ee \widetilde{T}_y(x,y,\xi,t)  \right)} \right]
d\xi  d\eta+  \nonumber   \\
&&\hskip-0.48in 
+ \frac{z+h}{h} \int_{-h}^0 \int_{-h}^{\eta}  \left[
\frac{\ee \Dd \widetilde{T} (x,y,\xi,t) 
- f_0 \widetilde{T}_x(x,y,\xi,t)}{\ee^2 +f^2} -
\right.   \nonumber      \\
&&\hskip-0.48in 
\left. - \frac{2 f_0 f}{(\ee^2 +f^2)^2} 
{\left( -f \widetilde{T}_x(x,y,\xi,t) 
+\ee \widetilde{T}_y(x,y,\xi,t)  \right)} \right]
d\xi  d\eta  \nonumber  \\
&&\hskip-0.48in
- \frac{ z(z+h)}{2} \left[
\frac{\ee \Dd T^* (x,y) - f_0 T^*_x(x,y)}{\ee^2 +f^2} 
+ \frac{2 f_0 f}{(\ee^2 +f^2)^2} 
{\left( -f T^*_x(x,y) +\ee T^*_y(x,y)  \right)} \right].  
\nonumber      
\end{eqnarray}
\bigskip

We denote by $L^p(\Om)$ and $H^m(\Om)$
the standard $L^p-$spaces and Sobolev spaces, respectively.
Following the notations in \cite{CT01} and \cite{STW98}, we set
\begin{equation}
\mathcal{V} = 
\left\{ { 
     R \in C^{\infty}(\overline{\Om}): 
     \left. { 
            \left( \pp_z \widetilde{T}+ \frac{\aa}{K_v} \widetilde{T} 
            \right)} 
            \right|_{z=0}= 0; \quad 
            \left. \pp_z \widetilde{T}  
            \right|_{z=-h}= 0; \quad 
            \left. q(  \widetilde{T}) \cdot \vec{n} 
     \right|_{\Gg_s}= 0 }
\right\}.  \label{FS-1}
\end{equation}
 For every $R \in C^{\infty}(\overline{\Om})$, denote by 
\begin{equation}
\| R \|_{V_2} = \left( \|R\|_{H^1 (\Om)}^2 + 
\|\nabla R\|_{H^1 (\Om)}^2 \right)^{1/2},
\label{N-2}
\end{equation}
and
\begin{equation}
\| R \|_{V_4} = \left( \|R\|_{H^2 (\Om)}^2 + 
\|\Dd R \|_{H^2 (\Om)}^2 \right)^{1/2}.
\label{N-4}
\end{equation}
Also, we denote by 
\begin{equation}
V_2 = {\mbox{the closure of }} \mathcal{V}
{\mbox{ with respect to the topology induced by the norm }} 
\| \cdot \|_{V_2}.       \label{FS-2}
\end{equation}
and
\begin{equation}
V_4 = {\mbox{the closure of }} \mathcal{V}
{\mbox{ with respect to the topology induced by the norm }} 
\| \cdot \|_{V_4}.   \label{FS-3}
\end{equation}
It is clear that $V_2$ and $V_4$ are separable Hilbert spaces. 
Next we define the bilinear form $a: V_2 \times V_2 \rightarrow \mathbb{R}$ 
is defined by 
\begin{eqnarray}
&&\hskip-.05in
a (R_1,R_2) =  \aa \int_{z=0} \left[ R_1 \, R_2 
+ \frac{\mu}{K_v} \nabla R_1 \cdot \nabla R_2 \right] dx dy + 
 \label{FORM}       \\
&&\hskip-.05in 
+ \int_{\Om} \left[  K_h  \nabla R_1 \cdot \nabla R_2 +
K_v (R_1)_z  \, (R_2)_z +\la \nabla \cdot (H^T \nabla R_1) \, 
\nabla \cdot (H^T \nabla R_2)
+ \mu \nabla (R_1)_z \cdot \nabla (R_2)_z 
\right] dxdydz.     \nonumber      
\end{eqnarray}
We will denote by 
\begin{eqnarray}
&&\hskip-.8in
| R |^2  =\int_{\Om} |R|^2 \; dxdydz,     \\
&&\hskip-.8in
\| R \|^2  = a(R,R) = \aa \int_{z=0} \left[ |R|^2 
+ \frac{\mu}{K_v} |\nabla R|^2 \right] dx dy +    \label{V-NORM}  \\
&&\hskip-.5in 
+ \int_{\Om} \left[  K_h  |\nabla R|^2 +
K_v |\pp_z R|^2 +\la |\nabla \cdot(H^T \nabla R) |^2 + \mu | \nabla R_z|^2 
\right] dxdydz.        \nonumber
\end{eqnarray}

Now we give the definition of weak and strong solutions to the model. 
\begin{definition} \label{D-1}
\thinspace Let $S$ be any fixed positive time. \thinspace
A function $\widetilde{T}(x,y,z,t)$ 
is called a weak solution of {\em (\ref{TEQ-1})--(\ref{TEQ-3})} 
on $[0,S]$ if
\[
\widetilde{T}   \in C_w([0,S], L^2(\Om)) \cap L^2([0,S], V_2),   
\quad
\widetilde{T}_t   \in L^1 ([0,S], V^{\prime}_2),
\]
and 
\begin{eqnarray}
&&\hskip-0.48in
\int_{\Om} \widetilde{T}(t) \psi \; dxdydz - \int_{\Om} 
\widetilde{T}(t_0)\psi  \; dxdydz  
+ \int_{t_0}^t a(\widetilde{T}(s), \psi) \; ds+  
\label{WEAK}  \\ 
&&\hskip-0.38in
+\int_{t_0}^t \int_{\Om}  \left[ v
 \cdot \nabla \widetilde{T}(s) +
w \widetilde{T}_z(s) 
+   v  \cdot \nabla T^*  \right] \psi \; dxdydz \; ds =\int_{t_0}^t 
\int_{\Om} Q^* \psi\;dxdydz \; ds,
\nonumber
\end{eqnarray}
for every $\psi \in V_2.$ 
Here $C_w([0,S], L^2(\Om))$ is the functional space of all weakly
continuous functions from $[0,S]$ to $L^2(\Om)$.
Furthermore, $\widetilde{T}(x,y,z,t)$ is
a strong solution of {\em (\ref{TEQ-1})--(\ref{TEQ-3})} on $[0,S]$ if
it is a weak solution and
\begin{eqnarray*}
&&\hskip-1.0in
\widetilde{T} (x,y,z,t)  \in C([0,S],V_2) 
\cap L^2([0,S], V_4).
\end{eqnarray*}
\end{definition}

Next, we give some remarks  about the following  
boundary value problem: 
\begin{equation}
\left\{ 
\begin{array}{l}
\nabla \cdot q (R) - K_v R_{zz} = g(x,y,z) 
\qquad {\mbox{in}}\quad  \Omega \\
\left. {\left( \pp_z R+ \frac{\aa}{K_v}  R 
 \right)} \right|_{z=0}= 0; \quad \left. \pp_z R
\right|_{z=-h}= 0; \quad \left. \frac{\pp R}{\pp \vec{e}} 
\right|_{\Gg_s}= 0;
\quad \left. q(  R) \cdot \vec{n} \right|_{\Gg_s}= 0,
\end{array}
\right.
\label{ELL}
\end{equation}
where $q(R)$ is given by (\ref{QQ}).
By integrating by parts and  the boundary conditions (\ref{HB}), we have
\begin{eqnarray}
&&\hskip-.8in
\int_{\Om} R_1 \left( \nabla \cdot q(R_2) - K_v (R_2)_{zz} \right)  dx dydz 
= \int_{\Om} R_2 \left( 
\nabla \cdot q(R_1) + \pp_z q_3(R_1) \right) dx dydz = a(R_1, R_2),
\label{SYMM}
\end{eqnarray}
for all $R_1, R_2 \in \widetilde{\mathcal{V}}.$   
Notice that 
\[
\vec{e} \cdot \vec{n} =\ee \not= 0.
\]
Namely, the vector $\vec{e}$ is not tangent to $\pp M.$ 
Using  the classical regularity results in smooth domains, 
$M \subset \mathbb{R}^2$, 
for the Laplacian operator with
oblique boundary condition, we have
\begin{proposition} \label{T-T}
Suppose that $R(\cdot, \cdot, z)$ satisfies the
boundary condition $\left. \frac{\pp R}{\pp \vec{e}} \right|_{\Gg_s} =0$
and that
$\Dd R(\cdot, \cdot, z)$ is in $L^2 (M)$ for every
fixed $z$. Then there exists a constant $C>0$ (independent of $z$) 
such that 
\begin{equation}
\|R(\cdot, \cdot, z) \|^2_{H^2(M)}  \leq C  
\left[ 
\|\Dd R (\cdot, \cdot, z) \|^2_{L^2(M )}
+\|R(\cdot, \cdot, z) \|^2_{L^2(M )}
\right].
\label{TT}
\end{equation}
Consequently, there is a constant $C_0>0$ such that
\begin{equation}
\frac{1}{C_0} \|R \|^2  \leq \| R \|_{V_2}^2 \leq C_0  \|R \|^2.
\label{POS}
\end{equation}
\end{proposition} 

\begin{proof}
As we mentioned earlier the proof of (\ref{TT}) is a result of classical 
regularity theory of elliptic equations. Notice that since  
$\displaystyle {\left. \frac{\pp R}{\pp \vec{e}} \right|_{\Gg_s} =0}$ 
(thanks to~(\ref{ELL})) we also have 
$\displaystyle{\left. \frac{\pp R_z}{\pp \vec{e}} \right|_{\Gg_s}=0}$. 
Now applying
(\ref{TT}) to $\partial_z R$ instead of $R$ with the corresponding boundary condition
$\displaystyle {\left. \frac{\pp R_z}{\pp \vec{e}} \right|_{\Gg_s}=0}$
 implies (\ref{POS}).
\end{proof}
Using (\ref{POS}),  Lax--Milgram
Theorem shows that there is a  unique solution $R \in V_2$ for
the boundary value problem (\ref{ELL}) satisfying
\begin{equation}
\|R\| \leq |g|.     \label{TEST-1}
\end{equation}
Moreover, using integration by parts in (\ref{ELL}) leads to
\begin{eqnarray*}
&&\hskip-.2in
|g|^2 = \left| \nabla \cdot q(R)-K_v R_{zz} \right|^2 
= \int_{\Om} \left[ \nabla \cdot q(R)
-K_v R_{zz} \right]^2 \; dxdydz   \\
&&\hskip-.2in
= \int_{\Om} \left[ |\nabla \cdot  q(R)|^2 + K_v^2 |R_{zz}|^2 
+ 2 K_v \nabla \cdot  q(R) R_{zz}  \right]  \; dxdydz   \\
&&\hskip-.08in
= \int_{\Om} \left[ |\nabla \cdot  q(R)|^2 + K_v^2 |R_{zz}|^2 
+ 2 K_v K_h |\nabla R_{z}|^2 + 2K_v \mu |\nabla R_{zz}|^2   
+ 2K_v \la |\nabla \cdot ( H^T \nabla R_{z})|^2 \right]  \; dxdydz + 
  \\
&&\hskip-.08in
+ 2 \aa \int_{z=0} \left[ \; K_h |\nabla R|^2 +
2\la |\nabla \cdot ( H^T \nabla R) |^2  \; \right]  dxdy.
\end{eqnarray*}
As a result, we have
\begin{equation}
| R_{zz}| \leq C |g |.     \label{TEST-2}
\end{equation}
Notice that the principle part of operator $\nabla \cdot q(R)-K_v R_{zz}$
is $-\la \Dd (-\Dd R + R) -K_v (-\Dd R + R)_{zz}$. Then, using a symmetry
argument in the $z$-direction and  the standard
regularity results for the Laplacian operator
(see,  for example, {\bf \cite{LADY}} p. 89), we get
\[
\|-\Dd R+ R\|_{H^2(\Om)} \leq C |g|.
\]
Therefore
\begin{equation}
\| R\|_{V_4} = \left(
\|\Dd R \|_{H^2(\Om)}^2  + \| R\|^2_{H^2(\Om)}
\right)^{1/2} \leq C |g|.
\label{EST}
\end{equation}
As a result of Proposition \ref{T-T}, (\ref{POS}) and
Rellich Lemma {\bf \cite{AR75}}, one can show that the operator
$\nabla \cdot  q(R)- K_v R_{zz}$
with domain $V_4$ is a positive self--adjoint operator with compact inverse.
Therefore, the space $L^2(\Om)$ possesses an orthonormal
basis $\{ \phi_k(x,y,z) \}_{k=1}^{\infty}$ of eigenfunctions of the
operator $\nabla \cdot  q(R)- K_v R_{zz}$, such that
\begin{equation}
\nabla \cdot  q(\phi_k)- K_v (\phi_k)_{zz}  = \ll_k \phi_k,   \label{ENG}
\end{equation}
where $0< \ll_1 \leq \ll_2 \leq \cdots,$ and
$\displaystyle{\lim_{k\rightarrow \infty}} \ll_k =\infty.$
Moreover, by stand results(cf. for example, \cite{C84}), we have
\begin{equation}
\frac{k}{C_1} \leq \frac{\ll_k}{\ll_1}.
\label{LL}
\end{equation}
We will denote by $H_m = \mbox{span} \{ \phi_1, \cdots, \phi_m \},$
and by $P_m: L^2(\Om) \rightarrow H_m$ the $L^2(\Om)$
orthogonal projection onto $H_m$.

\vskip0.1in

For convenience we recall the following classical
inequality  about the trace operator:
\begin{equation}
\| R \|_{L^2(\pp \Om)}^2 \leq C_2 \| R \|_{H^1(\Om)}^2,
\label{TRACE}
\end{equation}
and the following version of
Sobolev embedding and interpolation Theorems (cf. for example,
{\bf \cite{AR75}}):
\begin{eqnarray}
&&\hskip-.8in
\left\{ \begin{array}{l}
\| h(x,y) \|_{L^4(M)}  \leq C_3  \| h(x,y) \|_{L^2(M)}^{1/2}
\| h(x,y) \|_{H^1(M)}^{1/2},
\qquad \quad  \forall \; h \in H^1(M),      \\
\| h(x,y) \|_{L^q(M)}  \leq C_q \| h(x,y) \|_{H^1(M)},
\qquad \quad  \forall \; h \in H^1(M), \mbox{and } q < \infty,     \\
\| h(x,y) \|_{L^{\infty}(M)}  \leq C_3  \| h(x,y) \|_{L^2(M)}^{1/2}
\| h(x,y) \|_{H^2(M)}^{1/2},
\qquad \quad  \forall \; h \in H^2(M),      \\
\| h(x,y) \|_{H^1(M)}  \leq C_3  \| h(x,y) \|_{L^2(M)}^{1/2}
\| h(x,y) \|_{H^2(M)}^{1/2},
\qquad \quad  \forall \; h \in H^2(M),      \\
\|\Dd  h(x,y) \|_{L^4(M)}  \leq C_3  \| h(x,y) \|_{L^2(M)}^{3/8}
\| h(x,y) \|_{H^4(M)}^{5/8},
\qquad \quad  \forall \; h \in H^2(M).
\end{array} \right.
\label{SIT-2}   \\
&&\hskip-.8in
\left\{ \begin{array}{l}
\| g(x,y,z) \|_{L^3(\Om)}  \leq C_4
| g(x,y,z) |^{1/2} \| g(x,y,z) \|_{H^1(\Om)}^{1/2},  \\
\| g(x,y,z) \|_{L^6(\Om)}  \leq C_4
\| g(x,y,z) \|_{H^1(\Om)},
\end{array} \right.
\quad \forall \; g \in H^1(\Om).   \label{SIT-3}
\end{eqnarray}
Also, the following integral version of Minkowsky inequality 
for the $L^p$ spaces $p\geq 1$. 
Let $\Om_1 \subset \mathbb{R}^{m_1}$ and
 $\Om_2 \subset \mathbb{R}^{m_2}$ be two measurable sets, where
$m_1$ and $m_2$ are two positive integers. Suppose that
$f(\xi,\eta)$ is measurable over $\Om_1 \times \Om_2$. Then,  
\begin{equation}
\left[ { \int_{\Om_1} \left( \int_{\Om_2} |f(\xi,\eta)| d\eta 
\right)^p d\xi } \right]^{1/p}
\leq \int_{\Om_2} \left( \int_{\Om_1} |f(\xi,\eta)|^p d\xi 
\right)^{1/p} d\eta.
\label{MKY}
\end{equation}
Hereafter, $C,$ which may depend on the domain $\Om$ and the constant 
parameters $\ee, f_0, \beta, \aa, K_h, K_v, \la, \mu$ 
in the system (\ref{TEQ-1})--(\ref{TEQ-3}),  
will denote a constant that may change from 
line to line.

\section{Global Existence, Uniqueness  and well--posedness 
of Weak Solutions}    \label{S-2}
Now we are ready to show the global existence and uniqueness of 
weak solutions to the system (\ref{TEQ-1})--(\ref{TEQ-3}).
\begin{theorem} \label{T-WEAK}
Suppose that $T^* \in H^2(M)$ and $Q \in L^2(\Om).$ 
Then for every $T_0(x,y,z) \in L^2(\Om),$ and 
$S>0,$ there is a unique weak solution $\widetilde{T}$ of the system 
{\em  (\ref{TEQ-1})--(\ref{TEQ-3})}. Moreover, 
$\widetilde{T}$ satisfies
\begin{equation}
\pp_t \widetilde{T} \in L^{\frac{4 r}{3 r-1}}
(0,S; V_2^{\prime}),  \quad \forall \; \; r > 1,\quad
|\widetilde{T}(\ss)|^2 +
\int_{0}^{\ss} \|\widetilde{T} (s)\|^2 \; ds  \leq  
 K_2 (S, Q, \widetilde{T}_0, T^*), \quad \forall \; \ss \in [0, S],        
\label{L-W-I}
\end{equation}
where $V_2^{\prime}$ is the dual space of $V_2$, 
$K_1 (S, Q, \widetilde{T}_0, T^*),$ 
and $K_2 (S, Q, \widetilde{T}_0, T^*)$  
are as specified in (\ref{K-1}) and (\ref{K-2}), respectively.
\end{theorem}
\begin{proof} 
First, let us prove the existence of the weak solution for system 
(\ref{TEQ-1})--(\ref{TEQ-3}). We will use a Galerkin like procedure,
based on the eigenfunctions $\{ \phi_k\}_{k=1}^{\infty}$, to show the 
existence. Let $m \in \mathbb{Z}^+$ be fixed, the Galerkin 
approximating system of order $m$ that we use for
(\ref{TEQ-1})--(\ref{TEQ-3}) reads:
\begin{eqnarray}
&&\hskip-.38in
\frac{\pp}{\pp t} \widetilde{T}_m + 
\nabla \cdot q (\widetilde{T}_m) -K_v (\widetilde{T}_m)_{zz}   
+ P_m\left[{ v \cdot \nabla \widetilde{T}_m + w  
\frac{\pp \widetilde{T}_m }{\pp z}  
+ (v \cdot \nabla) T^*  }\right] 
= P_m Q^*, 
 \label{GTEQ-1}  \\ 
&&\hskip-.38in
\widetilde{T}_m(x,y,z,0) = P_m [T_0 (x,y,z)-T^*(x,y)],
\label{GTEQ-2}
\end{eqnarray}
where $\widetilde{T}_m =  \sum_{k=1}^m a_k(t) \phi_k(x,y,z)$, and 
$v(\widetilde{T}_m) =(v_1,v_2), w=w(\widetilde{T}_m)$ given in terms of 
$\widetilde{T}_m$ by the formulas below, and
\begin{eqnarray}
&&\hskip-.8in
q (\widetilde{T}_m)= -K_h \nabla \widetilde{T}_m + 
\la H \nabla (\nabla \cdot (H^T\nabla \widetilde{T}_m)) 
+ \mu \nabla (\widetilde{T}_m)_{zz},   
\label{GTEQ-3}    \\
&&\hskip-.8in
v_1 =  \int_{-h}^z 
\frac{\ee (\widetilde{T}_m)_x(x,y,\xi,t) 
+f (\widetilde{T}_m)_y(x,y,\xi,t)}{\ee^2 +f^2} d\xi +
\left(z+\frac{h}{2}\right) 
\frac{\ee T^*_x(x,y) +f T^*_y(x,y)}{\ee^2 +f^2}-  
\label{V1_G}        \\
&&\hskip-0.48in
- \frac{1}{h} \int_{-h}^0 \int_{-h}^{\eta} 
\frac{\ee (\widetilde{T}_m)_x(x,y,\xi,t) 
+f (\widetilde{T}_m)_y(x,y,\xi,t)}{\ee^2 +f^2} d\xi d \eta,
\nonumber        \\
&&\hskip-0.8in
v_2 =  \int_{-h}^z 
\frac{-f (\widetilde{T}_m)_x(x,y,\xi,t) 
+\ee (\widetilde{T}_m)_y(x,y,\xi,t)}{\ee^2 +f^2} d\xi 
+\left(z+\frac{h}{2}\right) 
\frac{-f T^*_x(x,y) +\ee T^*_y(x,y)}{\ee^2 +f^2} - 
\label{V2_G}         \\
&&\hskip-0.48in
- \frac{1}{h} \int_{-h}^0 \int_{-h}^{\eta} 
\frac{-f  (\widetilde{T}_m)_x(x,y,\xi,t) 
+\ee (\widetilde{T}_m)_y(x,y,\xi,t)}{\ee^2 +f^2} d\xi d\eta,
\nonumber      \\
&&\hskip-0.8in
w =  - \int_{-h}^z \int_{-h}^{\eta} \left[
\frac{\ee \Dd \widetilde{T}_m (x,y,\xi,t) 
- f_0 (\widetilde{T}_m)_x(x,y,\xi,t)}{\ee^2 +f^2} +
\right. 
\label{W_G}          \\
&&\hskip-0.48in 
\left. + \frac{2 f_0 f}{(\ee^2 +f^2)^2} 
{\left( -f (\widetilde{T}_m)_x(x,y,\xi,t) 
+\ee (\widetilde{T}_m)_y(x,y,\xi,t)  \right)} \right]
d\xi  d\eta+  \nonumber   \\
&&\hskip-0.48in 
+ \frac{z+h}{h} \int_{-h}^0 \int_{-h}^{\eta}  \left[
\frac{\ee \Dd \widetilde{T}_m (x,y,\xi,t) 
- f_0 (\widetilde{T}_m)_x(x,y,\xi,t)}{\ee^2 +f^2} -
\right.   \nonumber      \\
&&\hskip-0.48in 
\left. - \frac{2 f_0 f}{(\ee^2 +f^2)^2} 
{\left( -f (\widetilde{T}_m)_x(x,y,\xi,t) 
+\ee (\widetilde{T}_m)_y(x,y,\xi,t)  \right)} \right]
d\xi  d\eta  \nonumber  \\
&&\hskip-0.48in
- \frac{ z(z+h)}{2} \left[
\frac{\ee \Dd T^* (x,y) - f_0 T^*_x(x,y)}{\ee^2 +f^2} 
+ \frac{2 f_0 f}{(\ee^2 +f^2)^2} 
{\left( -f T^*_x(x,y) +\ee T^*_y(x,y)  \right)} \right].
\end{eqnarray}
We stress again that $v$ and $w$ depend on $m$ since they are functions of 
$\widetilde{T}_m$. However, we will drop the explicit dependence $m$ to 
simplify the notation.
The  equation (\ref{GTEQ-1}) is an ODE system with 
the unknown $a_k(t), k=1,\cdots, m.$ 
Furthermore, it is easy to check that each term of equation 
(\ref{GTEQ-1}) is locally Lipschitz in $\widetilde{T}_m$. 
Therefore, there is a unique solution 
$a_k(t), k=1,\cdots, m,$ to the equation (\ref{GTEQ-1}) for a short
interval of time $[0, S^*)$. 
By taking the $L^2(\Om)$ inner product of equation 
(\ref{GTEQ-1}) with $\widetilde{T}_m$, we obtain
\begin{eqnarray}
&&\hskip-.28in
\frac{1}{2} \frac{ d |\widetilde{T}_m|^2 }{d t} 
+  \| \widetilde{T}_m \|^2
+\int_{\Om} (v\cdot\nabla T^* ) \,
\widetilde{T}_m \, dxdydz +     \label{PROD}   \\
&&\hskip-.18in
+ \int_{\Om} P_m \left[ { v \cdot \nabla \widetilde{T}_m 
 w \pp_z \widetilde{T}_m } \right] \widetilde{T}_m \,  dxdydz
=\int_{\Om} \widetilde{T}_m Q^*   dxdydz.    \nonumber
\end{eqnarray}
It is easy to show by integrating by parts and 
by using the boundary conditions 
(\ref{B-1})--(\ref{B-3}) that 
\begin{equation}
\int_{\Om} P_m \left[{ v \cdot \nabla \widetilde{T}_m +
w \pp_z \widetilde{T}_m}\right] \widetilde{T}_m \;  dxdydz =0.
\label{EST-1}
\end{equation}
Furthermore, by H\"{o}lder inequality we have
\begin{eqnarray*}
&&\hskip-.18in
\left| { \int_{\Om} (v \cdot \nabla T^*)   
\widetilde{T}_m   dxdydz } \right|  
\leq C \| v \|_{L^6(\Om)} \| \nabla T^* \|_{L^3 (\Om)}
\; |\widetilde{T}_m |. \label{EST-2}
\end{eqnarray*}
By (\ref{SIT-3}), (\ref{V1_G}) and (\ref{V2_G}), we have
\begin{equation}
\| v \|_{L^6(\Om)} \leq C  \| v \|_{H^1(\Om)} 
\leq  C \left[ { \left\| \nabla \widetilde{T}_m \right\|_{H^1(\Om)} 
+ \|T^* \|_{H^2(M)}} \right]
\leq C \left[ { \left\| \widetilde{T}_m \right\| 
+ \|T^* \|_{H^2(M)}} \right].  \label{V6}
\end{equation}
By (\ref{SIT-2}), we obtain
\[
\| \nabla T^* \|_{L^3 (M)} \leq C  \|T^* \|_{H^2(M)}.
\]
As a result of the above estimate, we get
\begin{eqnarray}
&&\hskip-.18in
\left| { \int_{\Om} (v \cdot \nabla T^* )  
\widetilde{T}_m   dxdydz } \right|     \nonumber   \\
&&\hskip-.18in 
\leq C \left[ { 
\left\| \widetilde{T}_m \right\| 
+ \|T^* \|_{H^2(M)} } \right] 
  \| T^* \|_{H^2 (M)} \; |\widetilde{T}_m |. 
\nonumber   \\
&&\hskip-.18in 
\leq \frac{1}{4} \left\| \widetilde{T}_m \right\|^2 + 
C \|T^* \|_{H^2(M)}^2 \left( 1+ |\widetilde{T}_m |^2 \right). 
\label{EST-3}  
\end{eqnarray}  
Applying Cauchy--Schwarz inequality and the definition of $Q^*$, we obtain
\begin{eqnarray}
&&\hskip-.18in   
\left| \int_{\Om} Q^* \widetilde{T}_m  
\; dxdydz \right|  \nonumber   \\
&&\hskip-.18in  
\leq  \left| \int_{\Om} Q \widetilde{T}_m  
\; dxdydz \right|+ \left| \int_{\Om} \nabla \cdot q(T^*) \, 
 \widetilde{T}_m  \; dxdydz \right|   \nonumber   \\
&&\hskip-.18in 
\leq |Q| \; |\widetilde{T}_m | 
+ C \|T^* \|_{H^2(M)} \|\widetilde{T}_m \|  \nonumber     \\
&&\hskip-.18in 
\leq \frac{1}{2} |Q|^2 + \frac{1}{2} |\widetilde{T}_m |^2 
+ C \|T^* \|^2_{H^2(M)} + \frac{1}{4} \|\widetilde{T}_m \|^2.     
\label{EST-4}
\end{eqnarray}
Therefore, from the above estimates (\ref{EST-1})--(\ref{EST-4})  
and (\ref{POS}), (\ref{PROD}) gives
\begin{eqnarray}
&&\hskip-.18in
\frac{ d |\widetilde{T}_m|^2 }{d t} 
+ \| \widetilde{T}_m\|^2 
\leq C \left[ |Q|^2 + \|T^* \|^2_{H^2(M)} \right]  
+ C \|T^* \|^2_{H^2(M)} |\widetilde{T}_m|^2.     \label{T-STAR}  
\end{eqnarray}
Thanks to Gronwall inequality, we conclude 
\begin{eqnarray}
&&\hskip-.3in
|\widetilde{T}_m(t)|^2 \leq \left[ { |T_0|^2 
+ C \left( |Q|^2 + \|T^* \|^2_{H^2(M)} \right)  
} \right] \; e^{ \textstyle{ C \| T^* \|^2_{H^2(M)} \; t }} 
 \label{EST-L-2}
\end{eqnarray}
when $0 \leq t < S^*$. But since the right hand side is bounded 
as $t$ goes to $S^*$, we conclude that $\widetilde{T}_m(t)$ must exist
globally, i.e., $S^* = +\infty.$   
Therefore, for any given $S>0$ and any $t \in [0,S]$, we have 
\begin{eqnarray}
&&\hskip-.3in
|\widetilde{T}_m(t)|^2 \leq K_1 (S, Q, \widetilde{T}_0, T^*),
\label{WL-2} 
\end{eqnarray}
where
\begin{equation}
K_1 (S, Q, \widetilde{T}_0, T^*)   
=\left[ { |T_0|^2 
+ C \left( |Q|^2 + \|T^* \|^2_{H^2(M)}  \right) } \right] 
e^{ \textstyle{ C \| T^* \|^2_{H^2(M)} \; S }}.
 \label{K-1} 
\end{equation}
By integrating (\ref{T-STAR}) with respect to $t$ over $[0,S]$, and 
by (\ref{WL-2}), we get
\begin{eqnarray}  
&&\hskip-.3in
\int_{0}^S \|\widetilde{T}_m \|^2 \; ds  \leq   
 K_2 (S, Q, \widetilde{T}_0, T^*),        \label{WLL-2}
\end{eqnarray}
where
\begin{equation}
 K_2 (S, Q, \widetilde{T}_0, T^*) 
= |T_0|^2 + 
 C \left[ |Q|^2 + \|T^* \|^2_{H^2(M)} 
+ \|T^* \|^2_{H^2(M)} \, K_1 (S, Q, |T_0|, T^*)  
\right] \; S, \label{K-2} 
\end{equation}
and $K_1 (S, Q, \widetilde{T}_0, T^*)$ is as in (\ref{K-1}).
Notice that the estimate (\ref{WL-2}) is 
unbounded in time (i.e., as $S\rightarrow \infty$), 
but it is uniformly bounded in $m$. 
However, in Section \ref{S-A} we will present a sharper estimate which
is asymptotically bounded in time.
As a result of all the above we have  $\widetilde{T}_m$ exists 
globally in time and is uniformly bounded, in $m$, in 
$L^{\infty} ([0,S]; L^2(\Om))$ and $L^2 ([0,S]; V_2)$ norms.

Next, let us show that $\pp_t \widetilde{T}_m$ is uniformly bounded, 
in $m$, in the  $L^{\frac{4 r}{3 r-1}} ([0,S]; V_2^{\prime})$ norm for
every $r>1$.   {F}rom  (\ref{GTEQ-1}), we have, 
for every $\psi \in V_2$ 
\[
\ang{\frac{\pp}{\pp t} \widetilde{T}_m, \psi} = 
\ang{ P_m Q^*  - K_v (\widetilde{T}_m)_{zz} + \nabla \cdot 
q(\widetilde{T}_m)    
- P_m\left[{ v \cdot \nabla \widetilde{T}_m  
+ w \frac{\pp \widetilde{T}_m}{\pp z}  + (v\cdot\nabla ) T^* }\right], 
\psi}.
\]
Here, $\ang{\cdot, \cdot}$ is the dual action of 
$V_2^{\prime}$.  It is clear that
\begin{eqnarray}
&&\hskip-.5in
\left| \ang{P_m Q^*, \psi}\right| 
\leq \left| \ang{P_m Q , \psi}\right| + 
\left| \ang{\nabla \cdot q( T^*), \psi}\right|     \\ 
&&\hskip-.5in
\leq \|Q\| |\psi|
+ C \| T^*\|_{H^2(M)} \, \|\psi\|,  
\label{L2-1}  
\end{eqnarray}
and  by integration by parts we have       
\begin{eqnarray}
&&\hskip-.5in
\left| \ang{  - K_v (\widetilde{T}_m)_{zz} + \nabla \cdot 
q(\widetilde{T}_m),  \psi} \right| 
\leq C \|\widetilde{T}_m\| \; \|\psi\|. 
\label{L2-2}   
\end{eqnarray}
Next, let us get an estimate for
\begin{eqnarray*}
&&\hskip-0.25in
\left| \ang{ P_m\left[{ v \cdot \nabla ( \widetilde{T}_m+ T^* )} 
+w \frac{\pp \widetilde{T}_m}{\pp z} \right], \psi} \right|  \\ 
&&\hskip-0.18in 
= \left| \int_{\Om}  P_m\left[{ v \cdot \nabla ( \widetilde{T}_m+ T^* )} 
+w 
\frac{\pp \widetilde{T}_m}{\pp z} \right]\; \psi \; dxdydz \right|  \\
&&\hskip-0.18in
= \left| \int_{\Om}  \left[{ v \cdot \nabla ( \widetilde{T}_m+ T^* )}
+w
\frac{\pp \widetilde{T}_m}{\pp z} \right]\; \psi_m \; dxdydz \right|, 
\end{eqnarray*}
where $\psi_m = P_m \psi.$ Thus, by integration by parts, we obtain  
\begin{eqnarray}
&&\hskip-0.25in
\left| \ang{ P_m\left[{ v \cdot \nabla ( \widetilde{T}_m+ T^* )} 
+w
\frac{\pp \widetilde{T}_m}{\pp z} \right], \psi} \right|  
\label{ESTT}  \\ 
&&\hskip-0.18in
= \left| \int_{\Om}  \left[{ v \cdot \nabla \psi_m }
+w
\frac{\pp \psi_m }{\pp z} \right]\;  ( \widetilde{T}_m+
T^* ) \; dxdydz \right|.   \nonumber
\end{eqnarray}
Next, we estimate
\[
\left|{ \int_{\Om}  \left( v \cdot \nabla \psi_m \right) 
\left( \widetilde{T}_m +T^* \right)  \; dxdydz} \right|  
\leq  \| \nabla \psi_m \|_{L^6(\Om)} \; 
 \| v \|_{L^3(\Om)} \;
| \widetilde{T}_m +T^* |.
\]
Applying (\ref{N-2}), (\ref{SIT-3}) and Proposition \ref{T-T}, we have
\[
\| \nabla \psi_m \|_{L^6(\Om)} \leq C \| \nabla \psi_m \|_{H^1(\Om)}
\leq C \| \psi_m \|.
\]
By (\ref{SIT-3}) and (\ref{V6}), we reach
\[
\| v \|_{L^3(\Om)} \leq C \| v \|_{H^1(\Om)}^{1/2} \; | v |^{1/2}  
\leq C \left[ \| \widetilde{T}_m  \| + \| T^* \|_{H^2(M)} \right].
\]
Since
\[
|\widetilde{T}_m +T^* | \leq 
|\widetilde{T}_m| + h^{1/2}\|T^* \|_{L^2(M)},
\]
the above inequalities imply
\begin{eqnarray}
&&\hskip-0.8in 
\left|{ \int_{\Om}  \left( v \cdot \nabla \psi_m \right) 
\left( \widetilde{T}_m +T^* \right) 
}\right|   \label{DW-3}  \\ 
&&\hskip-0.8in  
\leq C \left[ \| \widetilde{T}_m \| + \| T^* \|_{H^2(M)} \right]\;
\left[  | \widetilde{T}_m | + \| T^* \|_{L^2(M)}
\right] \; \| \psi_m \|.     \nonumber
\end{eqnarray}
For the other term in (\ref{ESTT}), we use (\ref{WWWW}) to get
\begin{eqnarray*}
&&\hskip-0.25in 
\left|{  \int_{\Om} 
     w  \pp_z \psi_m \; \left(\widetilde{T}_m +T^* \right) \; dxdydz 
}\right|   \\
&&\hskip-0.25in 
= \left|{  \int_{\Om} 
  \left( \int_{-h}^z \nabla \cdot v \; d\xi \right)  
\pp_z \psi_m \; \left(\widetilde{T}_m +T^* \right) \; dxdydz 
}\right|   \\
&&\hskip-.25in 
\leq  \int_{\Om}   
\left( \int_{-h}^0 \left| \nabla \cdot v(x,y,\xi,t)\right| 
\,d\xi \right) 
\left| \pp_z \psi_m (x,y,z)\right|\; 
\left|\widetilde{T}_m(x,y,z,t) +T^*(x,y) \right| 
\;   dx dy dz  \\
&&\hskip-.25in   
= \int_M \left[ {  \int_{-h}^0 \left| 
   \nabla \cdot   v(x,y,\xi,t)\right| \,d\xi \;  
\int_{-h}^0  \left| \pp_z \psi_m (x,y,z)\right| \; 
 \left| \widetilde{T}_m (x,y,z,t) + T^* (x,y)\right|  \; dz } \right] \; dx dy.
\end{eqnarray*}
Using Cauchy--Schwarz inequality, we obtain
\begin{eqnarray*}
&&\hskip-0.25in 
 \int_{-h}^0 \left| \pp_z \psi_m (x,y,z)
\right| \; \left| 
\widetilde{T}_m (x,y,z,t) + T^* (x,y)\right| dz   \\
&&\hskip-0.18in 
 \leq  \left( { \int_{-h}^0 \left| \pp_z \psi_m (x,y,z)
\right|^2 dz }\right)^{\frac{1}{2}} 
\; \left( { \int_{-h}^0  \left| \widetilde{T}_m (x,y,z,t)+
T^* (x,y)\right|^2 dz    }\right)^{\frac{1}{2}}.
\end{eqnarray*}
Thus, 
\begin{eqnarray}
&&\hskip-.2in  
\left|{  \int_{\Om} 
     w  \pp_z \psi_m \; \left(\widetilde{T}_m +T^* \right) \; dxdydz 
}\right|     \nonumber  \\
&&\hskip-0.2in 
\leq \int_M    \left( \int_{-h}^0 |  \nabla \cdot  
v(x,y,\xi,t)| \,d\xi \right) 
\left( { \int_{-h}^0 | \pp_z \psi_m (x,y,z)
|^2 dz }\right)^{\frac{1}{2}} \left( { \int_{-h}^0  
| \widetilde{T}_m (x,y,z,t)+T^*(x,y) |^2 dz
    }\right)^{\frac{1}{2}}  \, dxdy 
\nonumber  \\
&&\hskip-0.2in 
\leq \left( { 
\int_M \left( { 
\int_{-h}^0  \left|  \nabla \cdot v (x,y,\xi,t) \right| d\xi 
}\right)^2 dx dy 
} \right)^{\frac{1}{2}}
 \left( { \int_M  \left( \int_{-h}^0  
\left| \pp_z \psi_m (x,y,z) \right|^2 dz \right)^{r^{\prime}} 
dx dy }\right)^{\frac{1}{2r^{\prime}}}  \times  \nonumber \\
&&\hskip-0.1in
\times  \left( { 
\int_M \left( { 
\int_{-h}^0  \left| \widetilde{T}_m (x,y,z,t) + T^* (x,y) \right|^2 dz 
}\right)^{r} dx dy 
}\right]^{\frac{1}{2 r}},   \label{TEMP} 
\end{eqnarray}
here we apply H\"{o}lder inequality and $1/r + 1/r^{\prime}=1,$ and 
$r>1.$ By using Minkowsky  inequality (\ref{MKY}),  we get
\begin{eqnarray*}
&&\hskip-0.5in  
\left[ { 
\int_M \left( { 
\int_{-h}^0  | \widetilde{T}_m (x,y,z,t) + T^* (x,y)|^2 dz 
}\right)^{r} dx dy 
}\right]^{\frac{1}{2r}}  \\
&&\hskip-0.45in  
\leq \left( \int_{-h}^0 \left( { \int_M  
|\widetilde{T}_m (x,y,z,t) + T^*(x,y)|^{2r} dx dy
}\right)^{\frac{1}{r}} dz  \right)^{\frac{1}{2}}.
\end{eqnarray*}
Thanks to (\ref{SIT-2}), for every fixed $z$ and $t$ we have
\begin{eqnarray*}
&&\hskip-0.25in  
\left( { \int_M  
\left|\widetilde{T}_m (x,y,z,t) + T^*(x,y)\right|^{2r^{\prime}} dx dy
}\right)^{\frac{1}{2 r}}  \leq 
C_r \left\| \widetilde{T}_m (z,t)+T^* 
\right\|_{L^2(M)}^{\frac{1}{2} + \frac{1}{2 r} } 
\left\|\widetilde{T}_m(z,t)+T^* 
\right\|_{H^2(M)}^{\frac{1}{2} - \frac{1}{2 r} }.
\end{eqnarray*}
As a result of the above and H\"{o}lder inequality, we obtain
\begin{eqnarray}
&&\hskip-0.25in 
\left( {  \int_{-h}^0 \left( { \int_M  
\left|\widetilde{T}_m (x,y,z,t) + T^*\right|^{2r} dx dy
}\right)^{\frac{1}{ r }} \; dz } \right)^{\frac{1}{2}}   
\nonumber    \\ 
&&\hskip-0.25in 
\leq C_r \left( {  \int_{-h}^0  \left\| \widetilde{T}_m (z,t)+
T^* \right\|_{L^2(M)}^2  \; dz } \right)^{\frac{r+1}{4 r} }
\left( {  \int_{-h}^0  \left\| \widetilde{T}_m (z,t)+
T^* \right\|_{H^2(M)}^2  \; dz } 
\right)^{\frac{r+1}{2 r} }   \nonumber  \\ 
&&\hskip-0.18in 
 \leq C_r \left| \widetilde{T}_m + T^* 
\right|^{\frac{r+1}{4 r} }
\left\| \widetilde{T}_m  + T^* 
\right\|_{V_2}^{\frac{r-1}{2 r} }.   \label{WWW}
\end{eqnarray}
Similarly, by using Minkowsky inequality (\ref{MKY}) and 
 (\ref{SIT-2}), we get
\begin{eqnarray*}
&&\hskip-0.25in 
\left( { \int_M  \left( \int_{-h}^0  
\left| \pp_z \psi_m (x,y,z) \right|^2 dz \right)^{r^{\prime}} 
dx dy }\right)^{\frac{1}{2r^{\prime}}}   \\
&&\hskip-0.18in 
 \leq \left( { \int_{-h}^0  \left( { \int_M \left| \pp_z \psi_m (x,y,z) 
\right|^{2r^{\prime}}
 \; dx dy }\right)^{1/r^{\prime}}\; dz} \right)^{1/2}   \\
&&\hskip-0.18in 
 \leq C_r \left( { \int_{-h}^0  \left\| \pp_z \psi_m  
    \right\|_{H^1(M)}^2  \; dz  }\right)^{1/2} \\
&&\hskip-0.18in 
 \leq C_r \| \psi_m \|.   
\end{eqnarray*}
Here we used Proposition \ref{T-T}. 
As a result of (\ref{TEMP}), (\ref{V1_G}), (\ref{V2_G}),  
(\ref{WWW}) and the above estimate, we have
\begin{eqnarray}
&&\hskip-0.25in 
\left|{ \int_{\Om} w \pp_z \psi_m\; (\widetilde{T}_m +T^*)\;    
dxdydz }\right| 
 \label{DW-4}     \\
&&\hskip-0.18in
\leq C_r \left( |\widetilde{T}_m| + \|T^*\|_{L^2(M)} 
\right)^{\frac{r+1}{2 r} } \;
\left( \|\widetilde{T}_m\| + \|T^*\|_{H^2(M)} 
\right)^{\frac{3 r-1}{2 r} }   \; 
\; \left\| \psi_m \right\|. \nonumber
\end{eqnarray}
By (\ref{DW-3}) and (\ref{DW-4}), we have
\begin{eqnarray*}
&&\hskip-.25in
\left| \ang{ P_m\left[{ v \cdot \nabla ( \widetilde{T}_m+ T^* )} 
+w \frac{\pp \widetilde{T}_m}{\pp z} \right], \psi} \right|     \\
&&\hskip-.18in
\leq C_r \left( |\widetilde{T}_m|  + \|T^*\|_{L^2(M)} 
\right)^{\frac{r+1}{2 r} } 
\left( \|\widetilde{T}_m\|  + \|T^*\|_{H^2(M)}
 \right)^{\frac{3 r-1}{2 r} }  \; \left\| \psi_m \right\|.    
\end{eqnarray*}
Since $\psi \in V_2$, then the Fourier series
\[
\sum_{k=1}^{\infty}\left( \int_{\Om} \psi \phi_k \; dxdydz \right)  
\phi_k = \psi_m + \sum_{k=m+1}^{\infty}
\left( \int_{\Om} \psi \phi_k \; dxdydz \right) \phi_k
\]
converges to $\psi$ in $V_2$ (cf. {\bf \cite{LADY}} p. 64).  
As a result, we get, 
\[
\left\| \psi_m \right\|_{V_2} \leq C \left\| \psi
\right\|_{V_2}.
\]
Therefore, from the above and Proposition \ref{T-T} we have
\begin{eqnarray}
&&\hskip-.25in
\left| \ang{ P_m\left[{ v \cdot \nabla ( \widetilde{T}_m+ T^* )} 
+w \frac{\pp \widetilde{T}_m}{\pp z} \right], \psi} \right| 
\label{L2-3}     \\
&&\hskip-.18in
\leq C_r \left( |\widetilde{T}_m|  + \|T^*\|_{L^2(M)} 
\right)^{\frac{r+1}{2 r} } 
\left( \|\widetilde{T}_m\|  + \|T^*\|_{H^2(M)}
\right)^{\frac{3 r-1}{2 r} }  \; \left\| \psi \right\|,  
\nonumber    
\end{eqnarray}
for every $r >1.$
By the estimates (\ref{L2-1})--(\ref{L2-3}), (\ref{WL-2}) and 
(\ref{WLL-2}),  we have
\begin{eqnarray*}
&&\hskip-.25in
\left| \ang{ \pp_t \widetilde{T}_m, \psi } \right| 
\leq  C \left( 
|Q| + \| T^*\|_{H^2(M)} \right)  \|\psi\| 
+  C \|\widetilde{T}_m\| \; \|\psi\|  +  \\
&&\hskip-.18in
+ C_r \left( |\widetilde{T}_m| + 
\|T^*\|_{L^2(M)}  \right)^{\frac{r+1}{2 r} } 
\left( \|\widetilde{T}_m\| + \|T^*\|_{H^2(\Om)}  
\right)^{\frac{3 r-1}{2 r} } \| \psi\|.
\end{eqnarray*}
Thus, due to (\ref{WL-2}) and (\ref{WLL-2}), we have
\begin{equation}
\int_0^S  \|\pp_t \widetilde{T}_m (t) 
\|_{V_2^\prime}^{\frac{4 r}{3 r-1} } dt 
\leq K_3 (S, Q, \widetilde{T}_0, T^*, r),  
\label{DL2}
\end{equation}
where 
\begin{eqnarray}
&&\hskip-.5in
 K_3 (S, Q, \widetilde{T}_0, T^*, r) 
= C \left[\; 1+ |Q| + \|T^* \|_{H^2(M)} \right]^{\frac{4r}{3r-1}} \; 
S + C K_2(S, Q, |T_0|, T^*) +  
\label{K-3}   \\
&&\hskip-.45in
+C_r \left[  K_1 (S, Q, \widetilde{T}_0, T^*) 
+ \|T^*\|_{H^2(\Om)} \right]^{\frac{2(r+1)}{3r-1}}
\; \left[ K_2 (S, Q, \widetilde{T}_0, T^*)  + 
\|T^* \|_{H^2(M)}^2 \; S\right].   
\nonumber 
\end{eqnarray}
Therefore, $\pp_t \widetilde{T}_m$ is uniformly bounded, 
in $m$, in the  $L^{\frac{4 r}{3 r-1} } ([0,S]; V_2^{\prime})$ norm, 
for every $r>1$.
Thanks to (\ref{WL-2}), (\ref{WLL-2}) and (\ref{DL2}), 
one can apply the Aubin's compactness Theorem 
(cf., for example, {\bf \cite{CF88}, \cite{LL69},
\cite{TT79}})  and extract a subsequence $\{ \widetilde{T}_{m_j} \}$ of 
$\{ \widetilde{T}_m \}$; 
a subsequence $\{ v_{m_j} \}$ 
of $\{ v_m = v(\widetilde{T}_m ) \}$
and a subsequence $\{ {\pp_t \widetilde{T}}_{m_j} \}$ 
of $\{ {\pp_t \widetilde{T}}_m \}$;
which converge to 
$\widetilde{T} \in L^{\infty} ([0,S]; L^2(\Om)) 
\cap L^2 ([0,S]; V_2)$ 
and ${\pp_t \widetilde{T}} \in L^{\frac{4 r}{3 r-1} } 
([0,S]; V_2^{\prime})$, respectively,
 in the following sense:
\[
\left\{ {
\begin{array}{ll}
\displaystyle{\widetilde{T}_{m_j} \rightarrow \widetilde{T}} 
& \displaystyle{ \mbox{in } L^{\infty} ([0,S]; L^2(\Om));}  \\
\displaystyle{ \widetilde{T}_{m_j} \rightarrow \widetilde{T}} 
&\displaystyle{ \mbox{in } L^2 ([0,S]; H^1(\Om)) \quad \mbox{weakly};} \\
\displaystyle{ \pp_t \widetilde{T}_{m_j} \rightarrow 
\pp_t \widetilde{T}} 
&\displaystyle{\mbox{in } L^{\frac{4 r}{3 r-1} } ([0,S]; V_2^{\prime} ) 
\quad \mbox{weakly}.} 
\end{array} } \right.  
\]
Notice that since $\widetilde{T}_{m_j} \in \mathcal{V}$ , 
by integration by parts it is clear that 
\begin{eqnarray*}
&&\hskip-0.38in
\int_{\Om} \widetilde{T}_{m_j}(x,y,z,t) P_{m_j} \psi \, dxdydz 
- \int_{\Om} 
\widetilde{T}_{m_j}(x,y,z,t_0) P_{m_j} \psi  \, dxdydz 
+   a(\widetilde{T}_{m_j}, P_{m_j} \psi)      \\
&&\hskip-0.38in 
+\int_{t_0}^t \int_{\Om}  \left[ { 
\left( v_{m_j} 
 \cdot \nabla \widetilde{T}_{m_j} \right) } \; P_{m_j}\psi 
+w_{m_j}  \pp_z \widetilde{T}_{m_j} \; P_{m_j}\psi
+   {\left( v_{m_j}  \cdot \nabla T^* \right)} P_{m_j}\psi \right]  \, 
dxdydz  \\
&&\hskip-0.38in 
=\int_{\Om} Q^* P_{m_j}\psi  \, dxdydz,    
\end{eqnarray*}
for every $\psi \in C^{\infty}([0,S]; V_2),$  
and for every $t$, $t_0\in [0,S]$.
By passing to the limit, one can show 
as in the case of Navier--Stokes 
equations (see, for example,  
{\bf \cite{CF88}, \cite{TT79}}) that $\widetilde{T}$ 
also satisfies (\ref{WEAK}).
In other words, $\widetilde{T}$ is a weak solution of the system 
(\ref{TEQ-1})--(\ref{TEQ-3}).

\vskip0.08in

Next, we show the uniqueness of the week solution.
Let $\widetilde{T}_1$ and $\widetilde{T}_2$ be two weak
solutions of the system (\ref{TEQ-1})--(\ref{TEQ-3}) with respect to initial 
values $\widetilde{T}_0^{\prime}$ and $\widetilde{T}^{\prime\prime}_0$, and let 
$( v^{\prime} = (v_1^{\prime}, v_2^{\prime}), w^{\prime})$ 
and $( v^{\prime\prime} = (v_1^{\prime\prime}, v_2^{\prime\prime}), 
w^{\prime\prime})$  be given by (\ref{V1_N})--(\ref{W_N}) with respect to 
$\widetilde{T}_1$ and  $\widetilde{T}_2$. Denote by 
$u =v^{\prime}-v^{\prime\prime}$, 
$u_3 =w^{\prime}-w^{\prime\prime}$
and $\chi =\widetilde{T}_1-\widetilde{T}_2.$
Unfortunately, since we only have 
\[
\pp_t \chi \in L^{\frac{4 r}{3 r-1}}
(0,S; V_2^{\prime}),  \qquad \forall \; \; r > 1,
\]
we are not able to apply the standard energy method and Lions Lemma 
(cf. {\bf \cite{TT79}} Lemma 1.2. p.260). 
However, let us instead consider
\[
\chi_m = P_m \chi,   
\qquad \mbox{and} \quad \widetilde{\chi}_m = \chi - P_m \chi
\]
where $P_m$ is the $L^2(\Om)$ orthogonal projection onto $H_m$. 
It is clear that $\chi_m$ satisfies the following equations:
\begin{eqnarray}
&&\hskip-.8in
\pp_t \chi_m  + \nabla \cdot q(\chi_m) 
- K_v (\chi_m)_{zz} 
+ P_m \left[ u \cdot \nabla \widetilde{T}_1   
+u_3 \pp_z \widetilde{T}_1  +
v^{\prime\prime} \cdot \nabla \chi 
+ w^{\prime\prime} \chi_z +    u  \cdot \nabla T^* \right] =0, 
 \label{UEQ-1}   \\
&&\hskip-.8in
\chi_m (x,y,z,0) =
P_m \left( 
\widetilde{T}^{\prime}_0 (x,y,z)- \widetilde{T}^{\prime\prime}_0 (x,y,z) 
\right), 
\label{UEQ-3} 
\end{eqnarray}
where $u=(u_1,u_2)$ and
\begin{eqnarray}
&&\hskip-0.8in
q (\chi_m)= -K_h \nabla \chi_m + 
\la H \nabla (\nabla \cdot (H^T \nabla \chi_m)) 
+ \mu \nabla (\chi_m)_{zz},   
\label{UEQ-4}    \\
&&\hskip-.8in
u_1 (x,y,z,t)=  \int_{-h}^z 
\frac{\ee \chi_x(x,y,\xi,t) 
+f \chi_y(x,y,\xi,t)}{\ee^2 +f^2} d\xi -  
\label{WV-1}        \\
&&\hskip-0.58in
- \frac{1}{h} \int_{-h}^0 \int_{-h}^z 
\frac{\ee \chi_x(x,y,\xi,t) 
+f \chi_y(x,y,\xi,t)}{\ee^2 +f^2} d\xi dz,
\nonumber        \\
&&\hskip-0.8in
u_2 (x,y,z,t)=  \int_{-h}^z 
\frac{-f \chi_x(x,y,\xi,t) 
+\ee \chi_y(x,y,\xi,t)}{\ee^2 +f^2} d\xi - 
\label{WV-2}         \\
&&\hskip-0.58in
- \frac{1}{h} \int_{-h}^0 \int_{-h}^z 
\frac{-f  \chi_x(x,y,\xi,t) 
+\ee \chi_y(x,y,\xi,t)}{\ee^2 +f^2} d\xi dz,
\nonumber      \\
&&\hskip-0.8in
u_3 (x,y,z,t)=  - \int_{-h}^z \int_{-h}^{\eta} \left[
\frac{\ee \Dd \chi (x,y,\xi,t) 
- f_0 \chi_x(x,y,\xi,t)}{\ee^2 +f^2} +
\right. 
\label{WW}          \\
&&\hskip-0.58in 
\left. + \frac{2 f_0 f}{(\ee^2 +f^2)^2} 
{\left( -f \chi_x(x,y,\xi,t) 
+\ee \chi_y(x,y,\xi,t)  \right)} \right]
d\xi  d\eta+  \nonumber   \\
&&\hskip-0.58in 
+ \frac{z+h}{h} \int_{-h}^0 \int_{-h}^z  \left[
\frac{\ee \Dd \chi (x,y,\xi,t) 
- f_0 \chi_x(x,y,\xi,t)}{\ee^2 +f^2} -
\right.   \nonumber      \\
&&\hskip-0.58in 
\left. - \frac{2 f_0 f}{(\ee^2 +f^2)^2} 
{\left( -f \chi_x(x,y,\xi,t) 
+\ee \chi_y(x,y,\xi,t)  \right)} \right]
d\xi  dz.  \nonumber 
\end{eqnarray}
By taking the $V_2^{\prime}$ dual action to 
equation (\ref{UEQ-1}) with $\chi_m$,  we obtain
\begin{eqnarray}
&&\hskip-.8in
\frac{1}{2} \frac{ d|\chi_m|^2}{dt}   + \| \chi_m\|^2 
= -  \ang{u \cdot \nabla \widetilde{T}_1   
+u_3 \pp_z \widetilde{T}_1  +
v^{\prime\prime} \cdot \nabla \chi 
+ w^{\prime\prime} \chi_z +    u  \cdot \nabla T^*, \chi_m} 
\nonumber  \\
&&\hskip-.8in
= - \int_{\Om} \left[ 
u \cdot \nabla \widetilde{T}_1   
+u_3 \pp_z \widetilde{T}_1  +
v^{\prime\prime} \cdot \nabla \chi 
+ w^{\prime\prime} \chi_z +    u  \cdot \nabla T^* \right] \; \chi_m
\; dxdydz.     \label{WTEQ-1}
\end{eqnarray}
Next we estimate the equation (\ref{WTEQ-1}) term by term.
\begin{enumerate}
\item 
\[
\left|{ \int_{\Om}  u \cdot \nabla T^* \chi_m \; dxdydz } \right|  
\leq  \| \nabla T^* \|_{L^6 (\Om)}  
\| u \|_{L^3 (\Om)} |\chi_m|.
\]
From (\ref{WV-1}) and (\ref{WV-2}), by applying (\ref{SIT-3}) and 
Proposition \ref{T-T}, we have
\[
\| u \|_{L^3 (\Om)} \leq C \| \nabla \chi \|_{L^3 (\Om)} \leq 
C \| \chi \|.
\]
Thus, by the above and Proposition \ref{T-T}, we obtain
\begin{equation}
\left|{ \int_{\Om}  u \cdot \nabla T^* \chi_m \; dxdydz } \right| 
\leq C \| T^* \|_{H^2(M)} \left( | \chi_m | \; \| \chi_m \| 
+ | \chi_m | \; \| \widetilde{\chi}_m \| \right).         \label{W-1}
\end{equation}

\item By (\ref{WV-1}) and (\ref{WV-2}) we get
\begin{eqnarray*}
&&\hskip-.28in 
\left|{ \int_{\Om}  u \cdot \nabla \widetilde{T}_1  
\chi_m } \; dxdydz \right|  
\leq  C \left|{ \int_{\Om}\left(\int_{-h}^0 |\nabla \chi (x,y,z,t)| dz\right) 
\; | \nabla \widetilde{T}_1| (x,y,z,t)\;  | \chi_m(x,y,z,t)|  } \; dxdydz \right|   \\
&&\hskip-.28in
=C \left|{ \int_{M}  \left(\left( \int_{-h}^0 |\nabla \chi(x,y,z,t)| dz\right) 
\; \int_{-h}^0 | \nabla \widetilde{T}_1(x,y,z,t)| \;   
| \chi_m(x,y,z,t)|  \, dz  \right) } \; dxdy \right|.
\end{eqnarray*}
Applying Cauchy--Schwartz inequality to the above, we obtain
\begin{eqnarray*}
&& \left|{ \int_{\Om}  u \cdot \nabla \widetilde{T}_1  
\chi_m } \; dxdydz \right| 
\leq C \int_M  \left[ {  \left( \int_{-h}^0 |\nabla \chi| \, dz \right)\;
\left( \int_{-h}^0 |\nabla  \widetilde{T}_1 |^2 \, dz 
\right)^{\frac{1}{2}}  
\; \left( \int_{-h}^0 |\chi_m|^2  \, dz  
\right)^{\frac{1}{2}} } \right] \; dxdy.
\end{eqnarray*}
By using H\"{o}lder inequality we reach
\begin{eqnarray*}
&&\hskip0.2in  \left|{ \int_{\Om}  u \cdot \nabla \widetilde{T}_1  
\chi_m } \; dxdydz \right| 
\leq C  \left[ { \int_M \left( \int_{-h}^0 |\nabla \chi| \, dz
    \right)^4 \; dxdy }\right]^{\frac{1}{4}}  \times  \\
&& \hskip0.3in  
\times 
 \left[ { \int_M  \left( \int_{-h}^0 |\nabla \widetilde{T}_1 |^2 \, dz 
\right) \; dxdy } \right]^{\frac{1}{2} }  \; 
\left[ { \int_M   \left( \int_{-h}^0 | \chi_m|^2  \, dz  \right)^2 \,
    dxdy }  \right]^{\frac{1}{4}}.
\end{eqnarray*}
By using Minkowsky  inequality (\ref{MKY})  we get
\[
\left[ { 
\int_M \left( { 
\int_{-h}^0  \left| \chi_m (x,y,z,t) \right|^2 dz 
}\right)^2 dx dy 
}\right]^{\frac{1}{2}} 
\leq \int_{-h}^0 \left( { \int_M \left|\chi_m (x,y,z,t) \right|^4 dxdy
}\right)^{\frac{1}{2}} dz.
\]
Thanks to (\ref{SIT-2}) for every fixed $z$ 
we have
\begin{eqnarray*}
&& 
\left({ \int_M   \left| \chi_m (x,y,z,t) \right|^4 dx dy 
}\right)^{\frac{1}{4}}  \leq 
C \left\| \chi_m (z,t) \right\|_{L^2(M)}^{\frac{1}{2}} 
\left\| \chi_m (z,t) \right\|_{H^1(M)}^{\frac{1}{2}}   \\
&&
\leq 
C \left\| \chi_m (z,t) \right\|_{L^2(M)}^{\frac{3}{4}} 
\left\| \chi_m (z,t) \right\|_{H^2(M)}^{\frac{1}{4}}.
\end{eqnarray*}
As for estimate (\ref{TT}), it is easy to obtain
\[
\left\| \chi_m (z,t) \right\|_{H^2(M)} 
\leq C \left( \| \Dd \chi_m (z,t) \|_{L^2(M)} 
+ \|\chi_m (z,t)\|_{L^2(M)} \right).
\]
As a result of the above we reach
\begin{eqnarray*}
&& 
\left({ \int_M   \left| \chi_m (x,y,z,t) \right|^4 dx dy 
}\right)^{\frac{1}{4}}  \leq 
C \left\| \chi_m (z,t) \right\|_{L^2(M)}^{\frac{3}{4}} 
\left\| \Dd \chi_m (z,t) \right\|_{L^2(M)}^{\frac{1}{4}}
+ C \|\chi_m (z,t)\|_{L^2(M)}^2.
\end{eqnarray*} 
Applying Young's inequality, we obtain
\begin{eqnarray*}
&& 
\int_{-h}^0 \left( {  \int_M  
 \left| \chi_m (x,y,z,t) \right|^4 dx dy 
}\right)^{\frac{1}{2}} dz   \\
&&
\leq 
C \int_{-h}^0  \left\| \chi_m (z,t) 
\right\|_{L^2(M)}^{\frac{3}{2}} 
\left\| \Dd \chi_m (z,t) \right\|_{L^2(M)}^{\frac{1}{2}} dz 
+|\chi_m (z,t)|^2    \\ 
&& 
\leq C \left( {
 \int_{-h}^0 \left\| \chi_m (z,t)\right\|_{L^2(M)}^2 \; dz  
} \right)^{\frac{3}{4}} \; \left( {
 \int_{-h}^0
\left\| \Dd \chi_m (z,t)\right\|^2_{L^2(M)}\; dz  
} \right)^{\frac{1}{4}}   \\
&&
\leq 
C \left| \chi_m (t) \right|^{\frac{3}{4}}
\left| \Dd \chi_m (t)\right|^{\frac{1}{4}} 
+C |\chi_m (t)|^2.
\end{eqnarray*}
Therefore, 
\begin{eqnarray}
&&  \left[ { 
\int_M \left( { 
\int_{-h}^0  \left| \chi_m (x,y,z,t) \right|^2 dz 
}\right)^2 dx dy 
}\right]^{\frac{1}{4}}  \leq C  
\left| \chi_m (t)\right|^{3/4}  
\left| \Dd \chi_m (t) \right|^{1/4}+C |\chi_m(t)|^2.
\label{NL-1}
\end{eqnarray}
As for obtaining (\ref{NL-1}),  by using (\ref{MKY}), 
(\ref{SIT-2}) and (\ref{TT}), and proposition \ref{T-T}, we have
\begin{eqnarray}
&& 
\left[ { 
\int_M \left( { 
\int_{-h}^0  \left| \nabla \chi (x,y,\xi,t) \right| d\xi 
}\right)^4 dx dy } \right]^{\frac{1}{4}} 
\leq C  \left| \chi  \right|^{1/4} \left| \Dd \chi \right|^{3/4}
+C |\chi|^{1/2} \|\chi\|^{1/2}.
\label{NL-2}
\end{eqnarray}
As a result of (\ref{NL-1}) and (\ref{NL-2}) we have
\begin{equation}
\left|{ \int_{\Om}  u \cdot \nabla \widetilde{T}_1  
\chi_m } \; dxdydz \right| 
\leq C \left\| \widetilde{T}_1 \right\| \; 
\left( | \chi_m| \; \| \chi_m \|  +  
\left| \widetilde{\chi}_m  \right|^{1/4} \| 
\widetilde{\chi}_m \|^{3/4} \; |\chi_m|^{3/4} \|\chi_m\|^{1/4}
\right).
 \label{W-2}
\end{equation}

\item Notice that from (\ref{WWWW}) we have
\[
u_3 = - \nabla \cdot \int_{-h}^z u(x,y,\xi,t) \; d\xi.
\]
By integrating by parts, we get
\begin{eqnarray*}
&&\hskip-.28in 
\left|{ \int_{\Om}  u_3 \, \pp_z\widetilde{T}_1 \, \chi_m } \; 
dxdydz \right|  
= \left|{ \int_{\Om}  \left( \int_{-h}^z u(x,y,\xi,t) \; d\xi \right)
 \, \cdot 
\left( \chi_m \; \nabla \pp_z \widetilde{T}_1 + 
\pp_z \widetilde{T}_1 \nabla \chi_m \right) } \; dxdydz \right|.
\end{eqnarray*}
By integration by parts we reach
\begin{eqnarray*}
&&\hskip-.28in 
 \left|{ \int_{\Om}  \int_{-h}^z u(x,y,\xi,t) \; d\xi  \, \cdot 
\nabla \chi_m  \pp_z \widetilde{T}_1  } \; dxdydz \right|  \\
&&\hskip-.2in 
=  \left|{ \int_{\Om} \left( u(x,y,z,t)  \cdot \nabla \chi_m 
+ \left( \int_{-h}^z u(x,y,\xi,t) \; d\xi\right)  
\, \cdot  \nabla \pp_z \chi_m \right)
 \widetilde{T}_1  } \; dxdydz \right|.
\end{eqnarray*}
Therefore, 
\begin{eqnarray*}
&&\hskip-.28in 
\left|{ \int_{\Om} u_3 \, \pp_z\widetilde{T}_1 \, \chi_m } \; 
dxdydz \right|
\leq  \left|{ \int_{\Om}  \left(\int_{-h}^0 |u| \; dz \right)\; 
|\nabla \pp_z \widetilde{T}_1 | \; |\chi_m|  } \; dxdydz \right| +  \\
&&\hskip-.2in 
 + \left|{ \int_{\Om}  |u| \; |\nabla \chi_m| \; |\widetilde{T}_1 | 
  } \; dxdydz \right| +
\left|{ \int_{\Om}  \left(\int_{-h}^0 |u| \; dz \right)\; 
|\widetilde{T}_1 | \; |\nabla \pp_z  \chi_m|  } \; dxdydz \right|.
\end{eqnarray*}
Following the steps of getting estimate (\ref{W-2}) we have
\begin{eqnarray}
&&\hskip-.28in 
\left|{ \int_{\Om}  \left(\int_{-h}^0 |u| \; dz \right)\; 
|\nabla \pp_z \widetilde{T}_1 | \; |\chi_m|  } \; dxdydz \right|  
 \label{W-3}         \\
&&\hskip-.2in  
= \left|{ \int_{M}  \left(\int_{-h}^0 |u| \; dz \right)\; 
\left(\int_{-h}^0 |\nabla \pp_z \widetilde{T}_1 | \; |\chi_m| \; dz \right) 
 } \; dxdy \right|   \nonumber   \\
&&\hskip-.2in  
\leq \left( \int_{M}  \left( \int_{-h}^0 |u| \; dz\right)^4 
\; dxdy \right)^{1/4}  
|\nabla \pp_z \widetilde{T}_1| \; 
\left( \int_{M}  \left( \int_{-h}^0 |\chi_m|^2 \; dz\right)^2 
\; dxdy \right)^{1/4}  \nonumber   \\
&&\hskip-.2in  
\leq C \left\| \widetilde{T}_1 \right\| 
\left( | \chi_m| \; \| \chi_m \|  +  
\left| \widetilde{\chi}_m  \right|^{1/4} \| 
\widetilde{\chi}_m \|^{3/4} \; |\chi_m|^{3/4} \|\chi_m\|^{1/4}
\right),
\nonumber
\end{eqnarray}
and
\begin{eqnarray}
&&\hskip-.28in 
\left|{ \int_{\Om}  |u| \; |\nabla \chi_m| \;
|\widetilde{T}_1 | } \; dxdydz \right|   \label{W-4}  \\
&&\hskip-.2in 
\leq C \left|{ \int_{M}  \left(\int_{-h}^0 |\nabla \chi| \; dz \right) \left( 
\int_{-h}^0 |\widetilde{T}_1 | \; |\nabla \chi_m| \; dz \right)
 } \; dxdy \right|   \nonumber   \\
&&\hskip-.2in 
\leq C \int_{M}  \left[ { 
\int_{-h}^0 | \nabla \chi| \; dz \max_{-h \leq z \leq 0} 
| \widetilde{T}_1(\cdot, \cdot, z, t)| \;   
\int_{-h}^0 | \nabla \chi_m| \; dz
} \right]  \; dxdy   \nonumber  \\ 
&&\hskip-.2in 
\leq C \int_{M}  \left[ { 
\int_{-h}^0 | \nabla \chi| \; dz   
\left( \int_{-h}^0 | \widetilde{T}_1|^2 \; dz  \right)^{1/4} 
\left( \int_{-h}^0 | (\widetilde{T}_1)_z|^2 \; dz  \right)^{1/4} 
\; \int_{-h}^0 | \nabla \chi_m| \; dz
} \right] \; dxdy   \nonumber  \\ 
&&\hskip-.2in
 \leq C \left| \widetilde{T}_1 \right|^{\frac{1}{2}}
\left\| \widetilde{T}_1 \right\|^{\frac{1}{2}}
\left[ | \chi_m |^{1/2} \;\| \chi_m \|^{3/2}
 + | \widetilde{\chi}_m  |^{1/4} \| \widetilde{\chi}_m \|^{3/4} 
\; |\chi_m|^{1/4} \|\chi_m\|^{3/4} \right].
  \nonumber 
\end{eqnarray}
and
\begin{eqnarray}
&&\hskip-.28in 
\left|{ \int_{\Om}  \int_{-h}^0 |u| \; dz \; |\nabla \pp_z \chi_m| \;
|\widetilde{T}_1|   } \; dxdydz \right|   \label{W-5}  \\
&&\hskip-.2in 
\leq C \int_{M} \left[ { \int_{-h}^0 | \nabla \chi| \; dz   
\left( \int_{-h}^0 | \nabla \pp_z \chi_m|^2 \; dz  \right)^{1/2} 
\left( \int_{-h}^0 |\widetilde{T}_1|^2 \; dz  \right)^{1/2} 
} \right] \; dxdy 
 \nonumber  \\ 
&&\hskip-.2in 
\leq C \left| \nabla \pp_z \chi_m \right|
 \left( \int_{M} \left( \int_{-h}^0 | \nabla \chi| \; dz  
\right)^4 \; dxdy    
\right)^{1/4}
\left( \int_M \left( \int_{-h}^0 | \widetilde{T}_1|^2 \; dz
  \right)^{2} \; dxdy  \right)^{1/4} 
 \nonumber  \\ 
&&\hskip-.2in
 \leq C \left| \widetilde{T}_1 \right|^{\frac{3}{4}}
\left\| \widetilde{T}_1 \right\|^{\frac{1}{4}}
\left| \chi \right|^{\frac{1}{4}} 
\left\| \chi \right\|^{\frac{3}{4}}\; \left\| \chi_m \right\|.   
 \nonumber 
\end{eqnarray}

\item By integrating by parts and (\ref{TEQ-3}), we reach
\[
\int_{\Om} \left[{  v^{\prime\prime} \cdot \nabla \chi +
w^{\prime\prime} \pp_z \chi  }\right] \chi_m \; dxdydz = 
\int_{\Om} \left[{  v^{\prime\prime} \cdot \nabla \widetilde{\chi}_m +
w^{\prime\prime} \pp_z \widetilde{\chi}_m  }\right] \chi_m \; dxdydz.
\]
Similar to get estimates (\ref{W-2})--(\ref{W-5}),  we get
\begin{eqnarray}
&&\hskip-.2in
\left| 
\int_{\Om} \left[{  v^{\prime\prime} \cdot \nabla \widetilde{\chi}_m +
w^{\prime\prime} \pp_z \widetilde{\chi}_m  }\right] \chi_m \; dxdydz
\right|    \label{W-6} \\
&&\hskip-.2in 
\leq C \left[ \| \widetilde{\chi}_m \| \; 
| \widetilde{T}_2  |^{1/4} \| 
\widetilde{T}_2 \|^{3/4} \; |\chi_m|^{3/4} \|\chi_m\|^{1/4}
+\| \widetilde{\chi}_m \| \; 
| \widetilde{T}_2  |^{1/4} \| 
\widetilde{T}_2 \|^{3/4} \; |\chi_m|^{3/4} \|\chi_m\|^{1/4}  
\right.    \nonumber \\
&&
\left. +| \widetilde{\chi}_m |^{\frac{1}{2}}
\| \widetilde{\chi}_m \|^{\frac{1}{2}}
| \widetilde{T}_2  |^{1/4} \| \widetilde{T}_2 \|^{3/4} 
\; |\chi_m|^{1/4} \|\chi_m\|^{3/4} 
+ | \widetilde{\chi}_m |^{\frac{3}{4}}
\| \widetilde{\chi}_m \|^{\frac{1}{4}}
| \widetilde{T}_2 |^{\frac{1}{4}} 
\| \widetilde{T}_2 \|^{\frac{3}{4}}\;
 \| \chi_m \| \right]  \nonumber 
\end{eqnarray}
\end{enumerate}

Therefore, by (\ref{W-1})--(\ref{W-6}), and Young's inequality, 
we have
\begin{eqnarray*}
&&\hskip-.2in
\frac{d|\chi_m|^2}{dt} 
\leq C  \left\|\widetilde{\chi}_m \right\|^2+        \\ 
&&\hskip-.1in 
+C \left[1+\| T^* \|_{H^2(M)}^2 +  \| \widetilde{T}_1 \|^2 
   + | \widetilde{T}_1 |^6 \, \| \widetilde{T}_1 \|^2
+ | \widetilde{T}_2  |^{2/3} \| \widetilde{T}_2 \|^2
\right] \; |\chi_m|^2 +   \\
&&\hskip-.1in 
+ C \left[ 1+ \| \widetilde{T}_1 \|^2 
+ | \widetilde{T}_1 |^6 \, \| \widetilde{T}_1 \|^2
+| \widetilde{T}_2  |^{2/3} \| \widetilde{T}_2 \|^2
\right] 
\; \left|\widetilde{\chi}_m \right|^2.
\end{eqnarray*}
Thanks to Gronwall inequality, we get
\begin{eqnarray*}
&&\hskip-.28in
|\chi_m(t)|^2
\leq \left[ |\chi_m(t_0)|^2  +  
C \int_{t_0}^t  \left\|\widetilde{\chi}_m (s)\right\|^2 \; ds
+ C \left[ \left( 1 + 
{\displaystyle{\max_{t_0 \leq s \leq t}}} | \widetilde{T}_1(s)|^6 
\right)
 \int_{t_0}^t  \| \widetilde{T}_1 (s) \|^2  \; ds  \right. \right.
+ \\
&&
\left. \left. + {\displaystyle{\max_{t_0 \leq s \leq t}}}  
| \widetilde{T}_2(s)  |^{2/3} 
\,\int_{t_0}^t  \| \widetilde{T}_2 (s) \|^2 \; ds 
\right] \; {\displaystyle{\max_{t_0 \leq s \leq t}}}
 \left|\widetilde{\chi}_m (s) \right|^2
\right]  \times \\
&& 
\times  e^{\textstyle{  
C \left[ 1+ {\displaystyle{\max_{t_0 \leq s \leq t}}} 
\left( | \widetilde{T}_1(s) |^6 + | \widetilde{T}_2(s) |^{2/3} 
\right) \left(  
\int_{t_0}^t \| \widetilde{T}_1 (s) \|^2 \; ds +
\int_{t_0}^t  \| \widetilde{T}_2 (s) \|^2  \; ds \right) \right]}} 
 \times  \\
&& 
\times e^{\textstyle{
C \left( 1+ \| T^* \|_{H^2(\Om)}^2  \right) (t-t_0) 
}}.  
\end{eqnarray*}
Notice that the right hand side is bounded, uniformly in $m$, 
for every fixed $t$. By passing the limit, using the Lebesgue
dominant convergence Theorem, we obtain
\begin{eqnarray}
&&\hskip-.28in
|\chi(t)|^2
\leq |\chi(t_0)|^2 \; e^{\textstyle{  
C \left[ 1+ {\displaystyle{\max_{t_0 \leq s \leq t}}} 
\left( | \widetilde{T}_1(s) |^6 + | \widetilde{T}_2(s) |^{2/3} 
\right) \left(  
\int_{t_0}^t \| \widetilde{T}_1 (s) \|^2 \; ds +
\int_{t_0}^t  \| \widetilde{T}_2 (s) \|^2  \; ds \right) \right]
}}  \times  \nonumber  \\
&& 
\times e^{\textstyle{
C \left( 1+ \| T^* \|_{H^2(\Om)}^2  \right) (t-t_0) }}.  \label{FIVE}
\end{eqnarray}
Recall that $ \widetilde{T}_1, \widetilde{T}_2 \in L^{\infty} ([0,S]; L^2(\Om)) 
\cap L^2([0,S]; V_2)$, therefore, (\ref{FIVE}) implies the continuous dependence 
of weak solutions on initial data as well as their uniqueness. 
\end{proof}

\begin{corollary}
The weak solution of the system {\em (\ref{TEQ-1})--(\ref{TEQ-3})} 
depends continuously on the initial data. That is, the problem is
well--posed.

\end{corollary}


\section{Global Existence, Uniqueness and Well--posedness 
of Strong Solutions}    \label{S-3}

In previous section we have proved the 
existence, uniqueness and well--posedness 
of the weak solution for the reformulated  system
(\ref{TEQ-1})--(\ref{TEQ-3}). In this section we show 
the global existence, uniqueness and well--posedness 
of strong solutions for the 
system (\ref{TEQ-1})--(\ref{TEQ-3}). 

\begin{theorem} \label{T-MAIN}
Suppose that $T^* \in H^4(M)$ and $Q \in L^2(\Om).$ 
Then for every $T_0(x,y,z)\in V_2$ and  $S>0,$ 
there is a unique strong solution  $\widetilde{T}$ 
of the system {\em  (\ref{TEQ-1})--(\ref{TEQ-3})}. 
Moreover, $\widetilde{T}$ satisfies
\begin{equation}
\|\widetilde{T} (t)\|^2 + \int_0^t \left|
\nabla \cdot q(\widetilde{T}) 
-K_v (\widetilde{T})_{zz} \right|^2 \; ds
\leq K_s (S, Q, T_0, T^*),         \label{L-S}
\end{equation}
where $0\leq t\leq S$ and 
$ K_s (S, Q, T_0, T^*)$ will be specified as in {\em (\ref{K-S})}.
\end{theorem}

\noindent
{\bf Remark:} The steps of the following proof are formal in the sense
that they can be made more rigorous by proving the corresponding
estimates first for the Galerkin approximation system 
based on the eigenfunctions of operator 
$\nabla \cdot q -K_v  (\widetilde{T})_{zz}$
with the boundary conditions (\ref{TEQ-2}). 
Then the estimates for the exact
solution can be established by passing to the limit in the Galerkin
procedure by using the appropriate ``Compactness Theorems''.

\begin{proof}
Let $\widetilde{T}$ be the weak solution with the initial 
datum $\widetilde{T}_0$. we will show that $\widetilde{T}$ 
is a strong solution if $T_0(x,y,z) \in V_2$.
First, we get {\it a priori} estimate for $\left| \widetilde{T}_z \right|$. 
Notice that
\[
- \int_{\Om} \widetilde{T}_{zz} \widetilde{T} \; dxdydz 
= \int_{\Om} \widetilde{T}_z^2 \; dxdydz 
+ \frac{\aa}{K_v} \int_{z=0} \widetilde{T}^2 \; dxdy,
\]
and
\begin{eqnarray*}
&&\hskip-.3in
- \int_{\Om} \widetilde{T}_{zz} \;  
\left[ \nabla \cdot q( \widetilde{T}) -K_v (\widetilde{T})_{zz} 
 \right] \; dxdydz       \\ 
&&\hskip-.28in
=\int_{\Om} \left[ \; 
K_h  | \nabla \widetilde{T}_z |^2 +
K_v |\widetilde{T}_{zz}|^2  +
\lambda  | \nabla \cdot (H^T \nabla   \widetilde{T}_z) |^2  
+ \mu  | \nabla \widetilde{T}_{zz} |^2  \; \right] \; dxdydz +  \\
&&\hskip-.2in
+ \frac{\aa \lambda}{K_v}  \int_{z=0} |\nabla \cdot (H^T \nabla 
\widetilde{T}) |^2 \; dxdy + \frac{\aa K_h}{K_v}  
\int_{z=0} |\nabla \widetilde{T} |^2 \; dxdy.
\end{eqnarray*}
By taking the inner product of equation (\ref{TEQ-1})  
with $\widetilde{T}_{zz}$ in $L^2(\Om)$, we reach
\begin{eqnarray*}
&&\hskip-.3in
\frac{1}{2} \frac{d  }{dt} \left( | \widetilde{T}_z|^2  
+ \frac{\aa}{K_v} \int_{z=0} \widetilde{T}^2 \; dxdy \right)
+ \frac{\aa }{K_v}  \int_{z=0} \left[ \lambda 
|\nabla \cdot (H^T \nabla  \widetilde{T}) |^2 
+ K_h |\nabla \widetilde{T} |^2 \right] \;
dxdy +   \\
&&\hskip-.2in 
 +  K_h  | \nabla \widetilde{T}_z |^2 + 
K_v  |\widetilde{T}_{zz}|^2 +
\lambda  | \nabla \cdot (H^T \nabla   \widetilde{T}_z) |^2 
 + \mu | \nabla \widetilde{T}_{zz} |^2            \\
&&\hskip-.3in  
=\int_{\Om} \left[ Q^* -  v \cdot \nabla \widetilde{T} 
- w \pp_z \widetilde{T} -  v \cdot \nabla T^* \right] 
\widetilde{T}_{zz} \; dxdydz.
\end{eqnarray*}
Let us consider the above equation term by term
\begin{enumerate}
\item It is clear that 
\begin{eqnarray}
&&\hskip-.3in 
\left| \int_{\Om} Q^*  \widetilde{T}_{zz} \; dxdydz
\right|= \left|\int_{\Om} \left(Q - \nabla \cdot q(T^*) \right) 
\widetilde{T}_{zz} \; dxdydz \right|     \label{SZ-1-1}         \\
&&\hskip-.3in 
\leq C \left( | Q | + \|T^* \|_{H^4(M) }
\right) | \widetilde{T}_{zz} |.
\nonumber
\end{eqnarray}

\item Integrating by parts we get
\begin{eqnarray*}
&&\hskip-.3in 
\int_{\Om} \left[ -  v \cdot \nabla \widetilde{T} 
- w \widetilde{T}_z -  v \cdot \nabla T^* \right] 
\widetilde{T}_{zz} \; dxdydz \\
&&\hskip-.3in 
= \int_{\Om} \left[ v_z \cdot \nabla \widetilde{T} 
+w_z \pp_z \widetilde{T} +  v_z \cdot \nabla T^* 
+ v \cdot \nabla \widetilde{T}_z +w \widetilde{T}_{zz} \right] 
\widetilde{T}_{z} \; dxdydz + \\
&&\hskip-.2in  
+ \frac{\aa}{K_v} \int_{z=0} \left[ v\cdot \nabla {\left( 
 \widetilde{T} +T^* \right) } \right] \widetilde{T} \; dxdy.
\end{eqnarray*}
Moreover,  by (\ref{EQ-3}) and boundary conditions 
(\ref{B-1})--(\ref{B-3}), we reach 
\begin{eqnarray} 
&&\hskip-.3in 
\int_{\Om} \left( 
v \cdot \nabla \widetilde{T}_z +w \widetilde{T}_{zz} \right) 
\widetilde{T}_{z} \; dxdydz =0 .    \label{SZ-1-2}
\end{eqnarray}
By using (\ref{V1_N})--(\ref{W_N}) and  boundary conditions
(\ref{TEQ-2}), we obtain
\begin{eqnarray} 
&&\hskip-.03in 
\int_{\Om} \left[ v_z \cdot \nabla \widetilde{T} 
+w_z  \widetilde{T}_z +  v_z \cdot \nabla T^* \right] 
\widetilde{T}_{z} \; dxdydz + 
\frac{\aa}{K_v} \int_{z=0} \left[ v\cdot \nabla {\left( 
 \widetilde{T} +T^* \right) } \right] \widetilde{T} \; dxdy   
\label{SZ-1-3}          \\
&&\hskip-.03in
= \int_{\Om} \frac{\ee}{\ee^2 +f^2} |\nabla \widetilde{T} + \nabla
T^*|^2 \;  \widetilde{T}_{z} \; dxdydz   
- \int_{\Om} \nabla \cdot v \; |\pp_z \widetilde{T}|^2 \; dxdydz +  
\nonumber      \\
&&\hskip.02in 
+\frac{\aa}{K_v} \int_{z=0} \left[ { \int_{-h}^0  
\frac{\ee}{\ee^2 +f^2} \left( 1+\frac{\xi}{h} \right) 
 |\nabla \widetilde{T}(x,y,\xi,t) + \nabla
T^*(x,y)|^2 \; d\xi  } \right] \; \widetilde{T} \; dxdy.
\nonumber   
\end{eqnarray}
Thanks to (\ref{SIT-3}) and (\ref{TT}) we have
\begin{eqnarray} 
&&\hskip-.03in 
\left| \int_{\Om} \frac{\ee}{\ee^2 +f^2} |\nabla \widetilde{T} + \nabla
T^*|^2 \;  \widetilde{T}_{z} \; dxdydz
\right| \leq C \left( \|\nabla \widetilde{T} \|_{L^4(\Om)}^2 +
\|\nabla T^* \|_{L^4(M)}^2 \right) 
\left| \widetilde{T}_z \right|  \label{SZ-1-4}   \\
&& 
\leq C \left( |\Dd \widetilde{T} |^2 +
\| T^* \|_{H^4(M)}^2 \right)\; 
\left| \widetilde{T}_z \right|  \nonumber   \\
&& 
\leq C \left( \| \widetilde{T} \|^2 +
\| T^* \|_{H^4(M)}^2 \right) 
| \widetilde{T}_z |.  \nonumber    
\end{eqnarray} 
Similarly, thanks to (\ref{SIT-2}) and (\ref{TT}) we have
\begin{eqnarray} 
&&\hskip-.03in 
\left| { \frac{\aa}{K_v} \int_{z=0} \left[ { \int_{-h}^0  
\frac{\ee}{\ee^2 +f^2} \left( 1+\frac{\xi}{h} \right) 
 |\nabla \widetilde{T}(x,y,\xi,t) + \nabla
T^*(x,y)|^2 \; d\xi  } \right] \; \widetilde{T} \; dxdy
} \right|           \label{SZ-1-5}    \\
&&
\leq C \left( \|\nabla \widetilde{T} \|_{L^4(\Om)}^2 +
\|\nabla T^* \|_{L^4(M)}^2 \right) 
\| \widetilde{T}( z=0) \|_{L^2(M)}    \nonumber   \\
&& 
\leq C \left( |\Dd \widetilde{T} |^2 +
\| T^* \|_{H^4(M)}^2 \right) 
\| \widetilde{T}(z=0) \|_{L^2(M)}  \nonumber   \\
&& 
\leq C \left( \| \widetilde{T} \|^2 +
\| T^* \|_{H^4(M)}^2 \right) 
\| \widetilde{T}(z=0) \|_{L^2(M)}.  \nonumber    
\end{eqnarray} 
Following the steps to get the estimate (\ref{W-4}) we obtain
\begin{eqnarray} 
&&\hskip-.03in 
\left| { \int_{\Om} \nabla \cdot v \; |\pp_z \widetilde{T}|^2 \;
    dxdydz } \right|           \label{SZ-1-6}    \\
&&
\leq C \left( |\Dd \widetilde{T} | +
\|\Dd T^* \|_{L^2(M)} \right) 
\| \widetilde{T}_z \|_{L^2(M)}\; 
 \| \nabla \widetilde{T}_z \|_{L^2(M)}    \nonumber   \\
&& 
\leq C \left( \| \widetilde{T} \| +
\| T^* \|_{H^4(M)} \right) 
\| \widetilde{T}_z \|_{L^2(M)} \; 
 \| \nabla \widetilde{T}_z \|_{L^2(M)}.  \nonumber    
\end{eqnarray} 
\end{enumerate}
Therefore, from (\ref{SZ-1-1})--(\ref{SZ-1-6}), we have
\begin{eqnarray*}
&&\hskip-.32in
\frac{1}{2} \frac{d  }{dt} \left(| \widetilde{T}_z|^2  
+ \frac{\aa}{K_v} \int_{z=0} \widetilde{T}^2 \; dxdy \right)
+ K_h  | \nabla \widetilde{T}_z |^2 + 
K_v  |\widetilde{T}_{zz}|^2 +
\lambda | \nabla \cdot (H^T \nabla    \widetilde{T}_z) |^2   
 + \mu | \nabla \widetilde{T}_{zz} |^2 
+   \\
&&\hskip-.25in 
 +\frac{\aa }{K_v}  \int_{z=0} \left[ \lambda 
|\nabla \cdot (H^T \nabla  \widetilde{T}) |^2 
+ K_h |\nabla \widetilde{T} |^2 \right] \; dxdy 
\\
&&\hskip-.32in  
\leq  C \left( | Q | + \|T^* \|_{H^4(M) }
\right) | \widetilde{T}_{zz} |
+ C \left( \| \widetilde{T} \|^2 +
\| T^* \|_{H^4(M)}^2 \right) 
| \widetilde{T}_z |+ \\
&&\hskip-.25in 
 + C \left( \| \widetilde{T} \|^2 +
\| T^* \|_{H^4(M)}^2 \right) 
\| \widetilde{T}(z=0) \|_{L^2(M)}
+
C \left( \| \widetilde{T} \| +
\| T^* \|_{H^4(M)} \right) 
\| \widetilde{T}_z \|_{L^2(M)} \; 
 \| \nabla \widetilde{T}_z \|_{L^2(M)}.
\end{eqnarray*}
By Cauchy--Schwartz inequality, we obtain
\begin{eqnarray}
&&\hskip-.53in
\frac{d  }{dt} \left( | \widetilde{T}_z|^2  
+ \frac{\aa}{K_v} \int_{z=0} \widetilde{T}^2 \; dxdy \right)
+ K_h  | \nabla \widetilde{T}_z |^2 + 
K_v  |\widetilde{T}_{zz}|^2 +
\lambda | \nabla \cdot (H^T \nabla    \widetilde{T}_z) |^2   
 + \mu | \nabla \widetilde{T}_{zz} |^2 
+   \label{Z-EST}        \\
&&\hskip-.45in 
 +\frac{\aa }{K_v}  \int_{z=0} \left[ \lambda 
|\nabla \cdot (H^T \nabla  \widetilde{T}) |^2 
+ K_h |\nabla \widetilde{T} |^2 \right] \; dxdy 
\nonumber     \\
&&\hskip-.53in  
\leq  C \left( | Q |^2 + \|T^* \|_{H^4(M) }^2
+ \| \widetilde{T} \|^2 \right) 
+  C \left( \| \widetilde{T} \|^2 +
\| T^* \|_{H^4(M)}^2 \right) \; \left(
\| \widetilde{T}(z=0) \|_{L^2(M)}^2
+ | \widetilde{T}_z |^2   \right).         \nonumber  
\end{eqnarray}
Again, by Gronwall inequality, we get
\begin{eqnarray}
&&\hskip-.3in
| \widetilde{T}_z(t)|^2  
+ \frac{\aa}{K_v} \int_{z=0} \widetilde{T}^2 (t)\; dxdy 
\leq K_z(t, Q, \widetilde{T}_0, T^*), 
  \label{KZ}     
\end{eqnarray}
where
\begin{eqnarray}
&&\hskip-.3in  
K_z(t, Q, \widetilde{T}_0, T^*) =  e^{ \textstyle{  
C \left(1+ \|T^*\|_{H^2(M)}^2 \right)\; t 
+ K_1(t,Q,  \widetilde{T}_0, T^*) }} \times   \label{K-Z}    \\
&&\hskip-.2in
\times  \left[ { \| T_0 \|^2+\|T^*\|_{H^2(M)}^2 
+K_1(t,Q, \widetilde{T}_0, T^*)+ 
\left( | Q |^2 + \|T^* \|_{H^4(M) }^2 \right)\; t }
\right],          \nonumber
\end{eqnarray}
and $K_1(t,Q, \widetilde{T}_0, T^*)$ is as in (\ref{K-1}).
Finally, let us show that
\begin{eqnarray*}
&&\hskip-1.0in
\widetilde{T} (x,y,z,t)  \in L^{\infty}([0,S],V_2) 
\cap L^2([0,S], V_4),
\end{eqnarray*}
By taking the inner product of equation (\ref{TEQ-1}) 
with $\nabla \cdot q( \widetilde{T}) -K_v \widetilde{T}_{zz}$ 
in $L^2(\Om)$ to get
\begin{eqnarray*}
&&\hskip-.3in
\frac{1}{2} \frac{d \| \widetilde{T}\|^2 }{dt} + \left|
\nabla \cdot q( \widetilde{T})-K_v\widetilde{T}_{zz}
 \right|^2   \\
&&\hskip-.3in 
=\int_{\Om} \left[ Q^* - v \cdot \nabla \widetilde{T} 
- w \pp_z \widetilde{T} - v \cdot \nabla T^* 
  \right] \left[  \nabla \cdot q( \widetilde{T}) -K_v
  \widetilde{T}_{zz}   \right] \; dxdydz.
\end{eqnarray*}
Note that
\begin{enumerate}

\item
\begin{eqnarray}
&&\hskip-.3in
\left| { \int_{\Om} Q^* \left[  \nabla \cdot q( \widetilde{T})
-K_v \widetilde{T}_{zz} \right] \; dxdydz } \right|  \label{SS-1}  \\
&&\hskip-.28in 
\leq \left|Q^*\right| \; |\nabla \cdot q( \widetilde{T}) 
-K_v \widetilde{T}_{zz} |   \nonumber   \\
&&\hskip-.28in 
 \leq \left[ |Q| + C \|T^*\|_{H^4(M)} \right] 
\; |\nabla \cdot q( \widetilde{T}) 
-K_v \widetilde{T}_{zz} |.
 \nonumber 
\end{eqnarray}

\item
\begin{eqnarray}
&&\hskip-.3in
\left| {  \int_{\Om} v \cdot \nabla T^* 
\left[  \nabla \cdot q( \widetilde{T}) -K_v \widetilde{T}_{zz}
\right] dxdydz  } \right|      \label{SS-2}  \\
&&\hskip-.28in 
\leq C |v| \| \nabla T^* \|_{L^{\infty} (M)} 
\left| \nabla \cdot q( \widetilde{T}) 
-K_v \widetilde{T}_{zz} \right|    \nonumber        \\
&&\hskip-.28in 
\leq C \|\widetilde{T} \| \; \| T^* \|_{H^4 (M)} 
\left| \nabla \cdot q( \widetilde{T})-K_v \widetilde{T}_{zz} 
\right|.        \nonumber 
\end{eqnarray}

\item~ Following similar steps to those which led to the 
estimate (\ref{W-4}) we have
\begin{eqnarray}
&&\hskip-.01in
\left| {  \int_{\Om} v \cdot \nabla \widetilde{T} 
\left[  \nabla \cdot q( \widetilde{T})
-K_v \widetilde{T}_{zz}   \right] \; dxdydz } \right|
\label{SS-3}   \\
&&
\leq C \left( 
| \nabla  \widetilde{T} | + \| \nabla T^* \|_{L^2 (M)}
\right) \left( 
\| \Dd  \widetilde{T} \|_{L^2 (\Om)} + \| T^* \|_{H^2 (M)}
\right)  
\left| \nabla \cdot q( \widetilde{T}) 
-K_v \widetilde{T}_{zz} \right|
\nonumber      \\
&&
\leq C \left(  \| T^* \|_{H^2 (M)}^2 + 
\|\widetilde{T} \|^2 \right) \;
\left|\nabla \cdot q( \widetilde{T}) 
-K_v \widetilde{T}_{zz} \right|. 
\nonumber 
\end{eqnarray}

\item~ Next, let us deal with the last term
\begin{eqnarray*}
&&\hskip-.3in
\left| {  \int_{\Om} w \pp_z \widetilde{T}
 \left[  \nabla \cdot q( \widetilde{T})
-K_v \widetilde{T}_{zz}   \right] \; dxdydz  } \right|      \\
&&\hskip-.28in 
\leq C \left\| w \right\|_{L^{\infty} (\Om)} \; 
| \pp_z \widetilde{T}| \; 
|  \nabla \cdot q( \widetilde{T}) -K_v \widetilde{T}_{zz} |.
\end{eqnarray*}
Notice that from  (\ref{W_N})
\begin{eqnarray*}
&&\hskip-.3in
\left\| w \right\|_{L^{\infty} (\Om)} 
\leq C \int_{-h}^0 \left[ 
\left\| \Dd  \widetilde{T} (z,t) \right\|_{L^{\infty} (M)} 
+ \left\| \Dd T^* \right\|_{L^{\infty} (M)} \right] \; dz 
\end{eqnarray*}
By using (\ref{SIT-2}) we obtain
\begin{eqnarray*}
&&\hskip-.3in
\| w \|_{L^{\infty} (\Om)}   
\leq C \left( {  \int_{-h}^0 
\|\Dd \widetilde{T} (z,t)\|_{L^2(M)}^{1/2} 
\| \Dd \widetilde{T} (z,t)\|_{H^2(M)}^{1/2}\; dz
+ \|\Dd T^* \|_{L^{\infty} (M)} }\right).  
\end{eqnarray*}
By Cauchy--Schwarz inequality, (\ref{EST}) and Proposition \ref{T-T} we get 
\begin{eqnarray*}
&&\hskip-.3in
\| w \|_{L^{\infty} (\Om)}   
\leq C \left( {  \|\Dd T^* \|_{L^{\infty} (M)} + 
\| \widetilde{T} \|^{1/2} \;
|  \nabla \cdot q( \widetilde{T}) -K_v \widetilde{T}_{zz} |^{1/2} } \right).
\end{eqnarray*}
As a result of the above we reach
\begin{eqnarray}
&&\hskip-.3in
\left| {  \int_{\Om} w \pp_z \widetilde{T}
 \left[  \nabla \cdot q( \widetilde{T}) -K_v \widetilde{T}_{zz}
  \right] \; dxdydz  } \right|  
\label{SS-4}  \\
&&\hskip-.28in 
\leq C \left[  \|T^*\|_{H^4 (M)} + \| \widetilde{T} \|^{1/2} \;
|  \nabla \cdot q( \widetilde{T}) -K_v \widetilde{T}_{zz}
 |^{1/2} \right] \times  \nonumber   \\
&&\hskip-.248in 
\times 
|\pp_z \widetilde{T}| \; 
|\nabla \cdot q( \widetilde{T})-K_v \widetilde{T}_{zz} |.         \nonumber
\end{eqnarray}
\end{enumerate}
Therefore, from the estimates (\ref{SS-1})--(\ref{SS-4}) 
we have
\begin{eqnarray*}
&&\hskip-.3in
\frac{1}{2} \frac{ d\| \widetilde{T}\|^2 }{dt} 
+ | \nabla \cdot q( \widetilde{T}) -K_v \widetilde{T}_{zz}|^2 \\
&&\hskip-.28in
\leq C \left|\nabla \cdot q( \widetilde{T}) 
-K_v \widetilde{T}_{zz} \right| \; 
\left[ 1+  |Q| + \|T^*\|_{H^4(M)}^2+ 
\|\widetilde{T} \|^2 + | \pp_z \widetilde{T}|^2 \right] +  \\
&&\hskip-.2in
 + C \| \widetilde{T} \|^{1/2} \;
| \pp_z \widetilde{T}| \; 
|\nabla \cdot q( \widetilde{T}) -K_v \widetilde{T}_{zz}|^{\frac{3}{2}}.  
\end{eqnarray*}
By using Young's inequality we get
\begin{eqnarray}
&&\hskip-.3in
\frac{d}{dt} \|\widetilde{T}\|^2 
+ \left| \nabla \cdot q( \widetilde{T})-K_v \widetilde{T}_{zz}
\right|^2   \label{V-EST}         \\ 
&&\hskip-.3in
\leq C \left[ 1+|Q|^2 + \|T^*\|_{H^4(M)}^4 +| \widetilde{T}_z|^4
\right] 
 + C \left[ 1  + 
\| \widetilde{T}\|^2 + | \widetilde{T}_z|^4 \right] 
\;  \| \widetilde{T}\|^2.      \nonumber
\end{eqnarray}
Again, by Gronwall inequality we conclude
\begin{eqnarray}
&&\hskip-.38in
\| \widetilde{T} (t)\|^2
+ \int_{0}^t  | \nabla \cdot q( \widetilde{T}) 
-K_v \widetilde{T}_{zz}
|^2 ds \leq  
K_s(t, Q, \widetilde{T}_0, T^*),    \label{KS}
\end{eqnarray}
where
\begin{eqnarray}
&&\hskip-.38in
K_s(S, Q, \widetilde{T}_0, T^*)= 
e^{ \textstyle{ C \left[ 1+ (K_z(t, Q, \widetilde{T}_0,  T^*))^2 
\right] t
+ K_1(t, Q, \widetilde{T}_0, T^*)   }} \times   \label{K-S}   \\
&&
 \times  \left[ \| \widetilde{T}_0 \|^2 + 
 |Q|^2 + \|T^*\|_{H^4(M)}^4 + K_z(t, Q, \widetilde{T}_0,  T^*)^2 \right],     
\nonumber
\end{eqnarray}
where $K_1(t, Q, \widetilde{T}_0, T^*)$ and 
$K_z(t, Q, \widetilde{T}_0, T^*)$ are as in 
(\ref{K-1}) and (\ref{K-Z}), respectively.
Since the strong solution is a weak solution, by Theorem \ref{T-WEAK},
the strong solution is unique.
\end{proof}

\section{Global Attractor}   \label{S-A}

In previous sections we have proved the existence and uniqueness 
of the weak and strong solution of the system
(\ref{TEQ-1})--(\ref{TEQ-3}). In this section we show the existence 
of the global attractor. Moreover, we give an upper bound, which are 
not necessarily optimal, for the
dimension of the global attractor. Denote by 
$\widetilde{T}(t) =S(t)\widetilde{T}_0$ 
the solution of the system (\ref{TEQ-1})--(\ref{TEQ-3}) with initial 
data $\widetilde{T}_0$. As a result of Theorem \ref{T-WEAK} and 
Theorem \ref{T-MAIN}, one can show that 
\[
\widetilde{T}(t) =S(t)\widetilde{T}_0 
\in L^2(\Om) \qquad \mbox{ for all} \quad  
\widetilde{T}_0 \in  L^2(\Om),  t \geq 0,
\]
and 
\[
\widetilde{T}(t) =S(t)\widetilde{T}_0 \in V_2 \qquad 
\mbox{ for all} \quad  \widetilde{T}_0 \in V_2, 
t \geq 0.
\] 
Since, in this section,  we only consider the long time behavior of 
solutions of the system (\ref{TEQ-1})--(\ref{TEQ-3}),  by 
Theorem \ref{T-WEAK} and Theorem \ref{T-MAIN}, we conclude that
$\widetilde{T}(t) \in L^{\infty}_{\mbox{loc}}((0,S], V_2)$ 
for every $\widetilde{T}(0) \in L^2(\Om).$
As a result, one can easily show that
\begin{equation}
\pp_t \widetilde{T} \in L^2_{\mbox{loc}}((0,S], V_2^{\prime}).
\label{SQR}
\end{equation}

\begin{theorem} \label{T-A}
Suppose that $Q \in L^2 (\Om)$ and $T^* \in H^4(M).$ 
Then, there is a global attractor $\mathcal{A} \subset  L^2(\Om)$ 
for the system {\em (\ref{TEQ-1})--(\ref{TEQ-3})}. 
Moreover, $\mathcal{A}$ has finite Hausdorff and fractal dimensions.
\end{theorem}

\begin{proof} 
First, let us show that there are absorbing balls in 
$L^2(\Om)$ and $V_2.$ 
Let $\widetilde{T}$ be the solution of the system 
(\ref{TEQ-1})--(\ref{TEQ-3}) with initial datum 
$\widetilde{T}_0 \in  L^2(\Om)$ and 
$|T_0| = \left| \widetilde{T}_0 +T^* \right| \leq
\rho.$ By Theorem \ref{T-WEAK} and Theorem \ref{T-MAIN}, there is a 
$t_0$ such that 
\[
|T(t_0)| = \left| \widetilde{T}(t_0) +T^* \right|
\leq 2 \rho.
\]
{F}rom now on we assume that $t\geq t_0.$ 
By taking the $H^{\prime}$ dual action to 
equation (\ref{TEQ-2}) with $\widetilde{T}$, we obtain
\begin{eqnarray*}
&&\hskip-.68in
\ang{\pp_t \widetilde{T}  + \nabla \cdot q(
  \widetilde{T})-K_v\widetilde{T}_{zz}, \widetilde{T} } 
+ \ang{ v \cdot \nabla \widetilde{T} +w
\pp_z \widetilde{T} + v  \cdot \nabla T^*, \widetilde{T} } =
\ang{Q^*, \widetilde{T} }. 
\end{eqnarray*}
Applying (\ref{SQR}) and Lions Lemma (cf. {\bf \cite{TT79}} 
Lemma 1.2. p.260), we reach 
\[
\ang{\pp_t \widetilde{T}, \widetilde{T} } 
= \frac{1}{2} \frac{d |\widetilde{T}|^2}{dt}.
\]
Moreover,  we have
\[
\ang{\nabla \cdot q( \widetilde{T})-K_v\widetilde{T}_{zz}, 
\widetilde{T} } = \|\widetilde{T}\|^2. 
\]
In addition, 
\begin{eqnarray*}
&&\hskip-.68in
 \ang{ v \cdot \nabla \widetilde{T}+w
\pp_z \widetilde{T} + v  \cdot \nabla T^*, \widetilde{T} }  
= \int_{\Om} \left[{ v \cdot \nabla \widetilde{T} 
+w \pp_z \widetilde{T} + v  \cdot \nabla T^* } \right]\; 
\widetilde{T} \; dxdydz;
\end{eqnarray*}
and
\[
\ang{Q^*, \widetilde{T}} = \int_{\Om} Q^* \, \widetilde{T} \; dxdydz,
\]
as long as the integrals on the right hand side make sense. 
Therefore, we have
\begin{eqnarray}
&&\hskip-.68in
\frac{1}{2} \frac{d |\widetilde{T}|^2}{dt}  + \| \widetilde{T} \|^2  
+ \int_{\Om} \left[{ v \cdot \nabla \widetilde{T}+w 
\pp_z \widetilde{T} + v  \cdot \nabla T^*} \right] 
\; \widetilde{T} \; dxdydz
= \int_{\Om} Q^* \, \widetilde{T} \; dxdydz.  \label{L-2-1}
\end{eqnarray}
By taking $\psi = T^*$  in the weak formulation (\ref{WEAK}), we get
\begin{eqnarray*}
&&\hskip-0.38in
\int_{\Om} \widetilde{T}(t) T^* \, dxdydz 
- \int_{\Om} \widetilde{T} (t_0) \, T^* \, dxdydz 
+  \int_{t_0}^t a(\widetilde{T}, T^*) +\\
&&\hskip-0.38in 
+\int_{t_0}^t \int_{\Om}  \left[ { v \cdot \nabla \widetilde{T} +w
 \pp_z \widetilde{T} + v  \cdot \nabla T^*   }\right]  \, T^* \; 
dxdydz =\int_{t_0}^t \int_{\Om} Q^* \, T^*   \; dxdydz.   
\end{eqnarray*}
It is equivalent to 
\begin{eqnarray}
&&\hskip-.68in
\frac{d }{dt}  \int_{\Om} \widetilde{T}(t) T^* \, dxdydz 
+ a(\widetilde{T}, T^*) +  
 \label{L-2-2}   \\
&&\hskip-.65in
+ \int_{\Om}  \left[ { v \cdot \nabla \widetilde{T}+w
 \pp_z \widetilde{T} + v  \cdot \nabla T^*  } \right]  \, T^* \; 
dxdydz   
=\int_{\Om} Q^* \, T^*   \; dxdydz.
 \nonumber  
\end{eqnarray}
By adding (\ref{L-2-1}) and (\ref{L-2-2}), we obtain
\begin{eqnarray*}
&&\hskip-.68in
\frac{1}{2} \frac{d  }{dt} \left( 
|\widetilde{T}|^2 + 2 \int_{\Om} \widetilde{T}(t)
  T^* \, dxdydz \right) + \| \widetilde{T} \|^2  
+ a(\widetilde{T}, T^*)  + \\
&&\hskip-.65in 
+ \int_{\Om} \left[{ v \cdot \nabla \widetilde{T} +w
\pp_z \widetilde{T} + v  \cdot \nabla T^* } \right] \; 
\left( \widetilde{T} +  T^* \right) \; dxdydz    \\
&&\hskip-.65in
= \int_{\Om} Q^* \, \left( \widetilde{T} +  T^* \right) \; dxdydz.  
\end{eqnarray*}
Notice that
\begin{enumerate}
\item~
\[
|\widetilde{T}|^2 + \int_{\Om} 2\widetilde{T} \, T^* \; dxdydz 
= |T|^2  - \int_{\Om} |T^*|^2  \, dxdydz;
\]
\item~
\begin{eqnarray*}
&&
 \int_{\Om} \left[{ v \cdot \nabla \widetilde{T} +w 
\pp_z \widetilde{T} + v  \cdot \nabla T^* } \right] \; 
\left( \widetilde{T} +  T^* \right) \; dxdydz       \\
&& = \int_{\Om}  \left[ v \cdot \nabla T +w \pp_z T \, 
\right] \; T \; dxdydz;
\end{eqnarray*}
\item~
\begin{eqnarray*}
&&
\int_{\Om} Q^* \, \left( \widetilde{T} +  T^* \right) \, dxdydz  \\
&&
= \int_{\Om}Q \, T \; dxdydz - a(\widetilde{T}, T^*) - 
\|T^*\|^2 + \aa \int_{z=0} \left[ \; T\, T^*  + 
\frac{\mu}{K_v} \nabla T\cdot \nabla T^* \; \right]\; dxdy.
\end{eqnarray*}
\end{enumerate}
Therefore, we get
\begin{eqnarray*}
&&\hskip-.68in
\frac{1}{2} \frac{d |T|^2}{dt} +
\|T\|^2 + \int_{\Om} \left[ v \cdot \nabla T +w  
\pp_z T \right] \, T \; dxdydz \\ 
&&\hskip-.65in 
= \int_{\Om} Q \, T \; dxdydz + \aa \int_{z=0} \left[ \; T\, T^*  + 
\frac{\mu}{K_v} \nabla T\cdot \nabla T^* \; \right]\; dxdy.
\end{eqnarray*}
By integration by parts and (\ref{B-3}), we obtain 
\begin{equation}
\int_{\Om} \left[ v \cdot \nabla T +w
 \pp_z T \right] \; dxdydz =0.
\label{E-1}
\end{equation}
By applying Cauchy--Schwarz inequality and (\ref{TRACE}), we reach 
\begin{eqnarray}
&&\hskip-.68in
\left| \int_{z=0} \left[ \; T\, T^*  + 
\frac{\mu}{K_v} \nabla T \cdot \nabla T^* \; \right]\; dxdy \right|   
\leq   \|T^*\|_{H^1(M)} \; \left[\; \| T\|_{L^2(M)}+ 
\frac{\mu}{K_v} \| \nabla T \|_{L^2(M)} \right]  \label{E-2}   \\
&&\hskip-.68in
\leq  C_2^{1/2} \|T^*\|_{H^1(M)} \; \left[\; \| T\|_{H^1(\Om)}+ 
\frac{\mu}{K_v} \| \nabla T \|_{H^1(\Om)} \right]   
\nonumber  \\
&&\hskip-.68in
\leq C_2^{1/2}  \left( 1+\frac{\mu}{K_v} \right) \|T^*\|_{H^1(M)} \|T\|. 
\nonumber
\end{eqnarray}
It is clear that 
\begin{equation}
\left| { \int_{\Om} Q \, T \; dxdydz } \right| 
\leq |Q|\; |T|.
\label{E-3}
\end{equation}
Therefore, by the above estimates (\ref{E-1})--(\ref{E-3}), we obtain
\begin{eqnarray*}
&&\hskip-.68in
\frac{d |T|^2}{dt} + \| T\|^2 
\leq C_2 \aa^2 \left( 1+\frac{\mu}{K_v} \right)^2
\|T^*\|^2_{H^1(M)} + 2 |Q| \; |T|.
\end{eqnarray*}
By Cauchy--Schwarz inequality and (\ref{POS}), we have
\begin{equation}
\frac{d |T|^2}{dt} +  \frac{1}{2} \| T\|^2 
\leq C_2 \aa^2 \left( 1+\frac{\mu}{K_v} \right)^2 
\|T^*\|^2_{H^1(M)} + 2 C_0 |Q|^2,    \label{EQU}
\end{equation}
where $C_0$ is as in (\ref{POS}). 
Thus, again, by (\ref{POS}), we obtain
\[
\frac{d |T|^2}{dt} + \frac{1}{ 2 C_0}  |T|^2
\leq C_2 \aa^2 \left( 1+\frac{\mu}{K_v} \right)^2
\|T^*\|_{H^1 (M)}^2 + 2 C_0 |Q|^2. 
\]
By Gronwall Lemma, we get
\begin{eqnarray*}
&&\hskip-.68in
|T(t)|^2 \leq |T(t_0)|^2 
e^{ \textstyle{ - \frac{1}{ 2 C_0}\, t }} 
+ 2C_0 C_2 \aa^2 \left( 1+\frac{\mu}{K_v} \right)^2
\|T^*\|_{H^1 (M)}^2 + 4 C_0^2 |Q|^2.
\end{eqnarray*}
As a result of the above, when $t$ is large enough such that 
\begin{eqnarray*}
&&\hskip-.68in
|T(t_0)|^2    
e^{\textstyle{ - \frac{1}{2C_0} \, t  }}  \leq 
2C_0 C_2 \aa^2 \left( 1+\frac{\mu}{K_v} \right)^2 
\|T^*\|_{H^1 (M)}^2 + 4C_0^2 |Q|^2,
\end{eqnarray*}
we have
\begin{eqnarray}
&&\hskip-.3in
 |T (t)|^2 \leq  \widetilde{R}_a(T^*, Q),   \label{A-LR-T}         
\end{eqnarray}
where
\begin{eqnarray}
&&\hskip-.3in
\widetilde{R}_a(T^*, Q)= 
4C_0 C_2 \aa^2 \left( 1+\frac{\mu}{K_v} \right)^2 
\|T^*\|_{H^1 (M)}^2 +  8C_0^2 |Q|^2.  
\label{R-A-T}         
\end{eqnarray}
In particular, 
\[
\limsup_{t \rightarrow \infty} |T(t)|^2 \leq 
2C_0 C_2 \aa^2 \left( 1+\frac{\mu}{K_v} \right)^2 
\|T^*\|_{H^1 (M)}^2 + 4C_0^2 |Q|^2.
\]
In other words, when $t$ is large enough we have
\begin{eqnarray}
&&\hskip-.3in
 |\widetilde{T} (t)|^2 = |T(t) -T^*|^2 \leq  R_a(T^*, Q),   
\label{A-LR}         
\end{eqnarray}
where
\begin{eqnarray}
&&\hskip-.3in
R_a(T^*, Q)=  2\widetilde{R}_a(T^*, Q) +2 \| T^*\|_{L^2(M)}^2.  
\label{R-A}         
\end{eqnarray}
where $\widetilde{R}_a(T^*, Q)$ is as  in (\ref{R-A-T}).
Therefore, there is an absorbing ball in $L^2(\Om)$ with radius
$R_a(T^*, Q)$ for system (\ref{TEQ-1})--(\ref{TEQ-3}).

\vskip0.1in

Next, we show that there is an absorbing ball in $V_2.$ 
First, notice that from (\ref{EQU}), we have
\begin{eqnarray*}
&&\hskip-.3in
\int_{t}^{t+r} \| T (s) \|^2 \; ds  \leq 
2  |T (t) | +  \left[
4C_0 C_2 \aa^2 \left( 1+\frac{\mu}{K_v} \right)^2 
\|T^*\|_{H^1 (M)}^2 +  8C_0^2 |Q|^2
\right] \; r 
\end{eqnarray*}
Therefore, by (\ref{A-LR}), when $t$ is large enough, we get
\begin{eqnarray}
&&\hskip-.3in
\int_{t}^{t+r} \|T(s) \|^2 \; ds  \leq  
 K_r (r, Q, T^*),     \label{KR}
\end{eqnarray}
where 
\begin{equation}
 K_r (r, Q, T^*) = 2 R_a(T^*, Q)+ \left[
4C_0 C_2 \aa^2 \left( 1+\frac{\mu}{K_v} \right)^2 
\|T^*\|_{H^1 (M)}^2 +  8C_0^2 |Q|^2
\right] \; r,  
\label{K-R} 
\end{equation}
and $R_a(T^*, Q)$ is as in (\ref{R-A}).

\vskip0.1in

{F}rom the proof of Theorem \ref{T-MAIN} we recall the inequality
(\ref{Z-EST}): 
\begin{eqnarray*}
&&\hskip-.3in
\frac{d }{dt} \left( | \widetilde{T}_z|^2  
+ \frac{\aa}{K_v} \int_{z=0} \widetilde{T}^2 \; dxdy \right)
 + K_h  | \nabla \widetilde{T}_z |^2 + 
K_v  |\widetilde{T}_{zz}|^2 +
\lambda | \nabla \cdot (H^T \nabla    \widetilde{T}_z) |^2   
 + \mu | \nabla \widetilde{T}_{zz} |^2 +        \\
&&\hskip-.2in 
+\frac{\aa }{K_v}  \int_{z=0} \left[ \lambda 
|\nabla \cdot (H^T \nabla  \widetilde{T}) |^2 
+ K_h |\nabla \widetilde{T} |^2 \right] \; dxdy     \\
&&\hskip-.3in  
\leq  C \left( | Q |^2 + \|T^* \|_{H^4(M) }^2
+ \| \widetilde{T} \|^2 \right) 
+  C \left( \| \widetilde{T} \|^2 +
\| T^* \|_{H^4(M)}^2 \right) \; \left(
\| \widetilde{T}(z=0) \|_{L^2(M)}^2
+ | \widetilde{T}_z |^2   \right).       
\end{eqnarray*}
Thus,
\begin{eqnarray*}
&&\hskip-.3in
\frac{d }{dt} \left( | \widetilde{T}_z|^2  
+ \frac{\aa}{K_v} \int_{z=0} \widetilde{T}^2 \; dxdy \right)    \\
&&\hskip-.3in  
\leq  C \left( | Q |^2 + \|T^* \|_{H^4(M) }^2
+ \| \widetilde{T} \|^2 \right) 
+  C \left( \| \widetilde{T} \|^2 +
\| T^* \|_{H^4(M)}^2 \right) \; \left(
\frac{\aa}{K_v} \| \widetilde{T}(z=0) \|_{L^2(M)}^2
+ | \widetilde{T}_z |^2   \right).       
\end{eqnarray*}
By applying the uniform Gronwall inequality (cf. for example, {\bf
  \cite{TT88}}, p. 89) and (\ref{KR}), we obtain, when $t$ is large enough
\begin{equation}
| \widetilde{T}_z(t)|^2 + \frac{\aa}{K_v} \| \widetilde{T}(z=0)
  \|_{L^2(M)}^2 
\leq R_z (r,T^*,Q),   \label{EST-Z}
\end{equation}
where $r>0$ is fixed and
\begin{eqnarray}
&&\hskip-.3in
R_z (r,T^*,Q) = 
C \left[ { 
\frac{K_r(r, T^*, Q)}{r} + \| T^*\|_{H^4 (M)}^2\; r +
|Q|^2\; r + K_r(r, T^*, Q) } \right] \times  \label{R-Z}   \\
&&\hskip-.25in
\times
e^{\displaystyle{ C \left[ K_r (r, T^*, Q) 
+ \|T^*\|_{H^4(M)}^2 \;r 
\right]  }}.   \nonumber  
\end{eqnarray}
Let us recall the inequality (\ref{V-EST}): 
\begin{eqnarray*}
&&\hskip-.18in
\frac{d}{dt} \|\widetilde{T}\|^2 
+ \left| \nabla \cdot q( \widetilde{T})-K_v \widetilde{T}_{zz}
\right|^2            \\ 
&&\hskip-.3in
\leq C \left[ 1+|Q|^2 + \|T^*\|_{H^4(M)}^4 +| \widetilde{T}_z|^4
\right] 
 + C \left[ 1  + 
\| \widetilde{T}\|^2 + | \widetilde{T}_z|^4 \right] 
\;  \| \widetilde{T}\|^2. 
\end{eqnarray*}
Thus, 
\begin{eqnarray*} 
&&\hskip-.3in
\frac{ d\| \widetilde{T} \|^2 }{dt} 
\leq C \left[ 1+|Q|^2 + \|T^*\|_{H^4(M)}^4 +| \widetilde{T}_z|^4
\right]  + C \left[ 1  + 
\| \widetilde{T}\|^2 + | \widetilde{T}_z|^4 \right] 
\;  \| \widetilde{T}\|^2
\end{eqnarray*}
Thanks the uniform Gronwall inequality (\ref{KR}) and 
(\ref{EST-Z}), we obtain, when $t$ is large enough
\begin{equation}
\| \widetilde{T}(t) \|
\leq R_v (r,T^*,Q),   \label{EST-V}
\end{equation}
where $r>0$ is fixed and
\begin{eqnarray}
&&\hskip-.3in
(R_v (r,T^*,Q))^2 = 
C \left[\frac{K_r(r,Q,T^*)}{r}+  r+|Q|^2 \;r
+ \|T^*\|_{H^4(M)}^4\;r + R_z(r,Q,T^*)^2\;r  \right]  
\times \label{R-V}   \\
&&\hskip-.25in 
\times 
e^{\displaystyle{  C \left[ r  + R_z(r,Q,T^*)^2 \; r +
K_r(r,Q,T^*)   \right]  }}.   \nonumber  
\end{eqnarray}
Therefore, we have shown that there is an absorbing ball $\mathcal{B}$
in $V_2$ with radius $R_v (r,T^*, Q)$. {F}rom the proofs of Theorems 
\ref{T-WEAK} and \ref{T-MAIN} 
we conclude that the operator $S(t)$ is a compact operator. 
Following the standard procedure (cf., for example, {\bf \cite{CF85}},
{\bf \cite{CF88}},
{\bf \cite{EFNT}, \cite{LADY91}, \cite{TT88}} for details), 
one can prove that there is a global 
attractor 
\[
\mathcal{A} =\cap_{t>0} S(t) \mathcal{B}  \subset V_2.
\]
Moreover, $\mathcal{A}$ is compact in $L^2(\Om)$ due to the compact 
embedding of $V_2$ in $L^2(\Om)$. 

\vskip0.1in

In addition to the compactness of the semi--group $S(t)$ one can 
show its differentiability on 
$\mathcal{A}$ with respect to the initial data. Therefore, one can 
use the trace formula (cf. {\bf \cite{CF85}, \cite{CF88}, 
\cite{TT88}}) to get an upper  bound for the dimension of the 
global attractor $\mathcal{A}$. Let $\widetilde{T}$ be a 
given solution of the system (\ref{TEQ-1})--(\ref{TEQ-3}) with 
$\widetilde{T} \in \mathcal{A}.$ Since it is on the global 
attractor $\mathcal{A}$, $\widetilde{T}$ is a strong solution 
to the system (\ref{TEQ-1})--(\ref{TEQ-3}).
It is clear that the first variation equations of the
 system (\ref{TEQ-1})--(\ref{TEQ-3}) around 
$\widetilde{T}$ read:
\begin{eqnarray}
&&\hskip-.68in
\pp_t \chi = F^{\prime} (\widetilde{T}) \chi,   \label{DIM-1}  \\ 
&&\hskip-.8in
\left. {\left( \pp_z \chi+ \frac{\aa}{K_v}  \chi
 \right)} \right|_{z=0}= 0; \quad \left. \pp_z \chi  
\right|_{z=-h}= 0; \quad \left. \frac{\pp \chi}{\pp \vec{e}} 
\right|_{\Gg_s}= 0;
\quad \left. q(  \chi) \cdot \vec{n} \right|_{\Gg_s}= 0,      
   \label{DIM-2} \\
&&\hskip-.68in
\chi (x,y,z,0) = \zeta,
 \label{DIM-3} 
\end{eqnarray}
where $\chi$ are the unknown perturbations about  
$\widetilde{T}$ with a given initial perturbation 
$\zeta \in L^2(\Om)$. Moreover, here 
\begin{eqnarray*}
&&\hskip-.68in
F^{\prime} (\widetilde{T}) \chi = - \left[
\nabla \cdot q(\chi)-K_v \chi_{zz} + u
\cdot \nabla (\widetilde{T}+T^*) + v \cdot \nabla \chi - 
 \widetilde{w}  \pp_z \widetilde{T}+ w \pp_z \chi 
\right],
\end{eqnarray*}
where $v=(v_1,v_2), w$ are as in (\ref{V1_N})--(\ref{W_N}), and 
\begin{eqnarray}
&&\hskip-0.8in
u_1 =  \int_{-h}^z 
\frac{\ee \chi_x(x,y,\xi,t) 
+f \chi_y(x,y,\xi,t)}{\ee^2 +f^2} d\xi -  
\label{DV-1}        \\
&&\hskip-0.58in
- \frac{1}{h} \int_{-h}^0 \int_{-h}^z 
\frac{\ee \chi_x(x,y,\xi,t) 
+f \chi_y(x,y,\xi,t)}{\ee^2 +f^2} d\xi dz,
\nonumber        \\
&&\hskip-0.8in
u_2 =  \int_{-h}^z 
\frac{-f \chi_x(x,y,\xi,t) 
+\ee \chi_y(x,y,\xi,t)}{\ee^2 +f^2} d\xi - 
\label{DV-2}         \\
&&\hskip-0.58in
- \frac{1}{h} \int_{-h}^0 \int_{-h}^z 
\frac{-f  \chi_x(x,y,\xi,t) 
+\ee \chi_y(x,y,\xi,t)}{\ee^2 +f^2} d\xi dz,
\nonumber      \\
&&\hskip-0.8in
\tilde{w} =  - \int_{-h}^z \int_{-h}^{\eta} \left[
\frac{\ee \Dd \chi (x,y,\xi,t) 
- f_0 \chi_x(x,y,\xi,t)}{\ee^2 +f^2} +
\right. 
\label{DW}          \\
&&\hskip-0.58in 
\left. + \frac{2 f_0 f}{(\ee^2 +f^2)^2} 
{\left( -f \chi_x(x,y,\xi,t) 
+\ee \chi_y(x,y,\xi,t)  \right)} \right]
d\xi  d\eta+  \nonumber   \\
&&\hskip-0.58in 
+ \frac{z+h}{h} \int_{-h}^0 \int_{-h}^z  \left[
\frac{\ee \Dd \chi (x,y,\xi,t) 
- f_0 \chi_x(x,y,\xi,t)}{\ee^2 +f^2} -
\right.   \nonumber      \\
&&\hskip-0.58in 
\left. - \frac{2 f_0 f}{(\ee^2 +f^2)^2} 
{\left( -f \chi_x(x,y,\xi,t) 
+\ee \chi_y(x,y,\xi,t)  \right)} \right],
d\xi  dz  \nonumber 
\end{eqnarray}
It is not difficult to show that the above, coupled  second order elliptic
 and  linear parabolic, system  has a unique
 solution $\chi(t)$. Moreover,  for $t>0$,
\begin{eqnarray*}
&&\hskip-.68in
\chi(t) \in V_2.
\end{eqnarray*}

For any positive integer $m$ we consider the volume element 
$\left|\chi_1(t)\wedge \cdots \wedge \chi_m(t)
\right|_{\wedge^{m} L^2(\Om)}$, we have the following trace 
formula (cf. {\bf \cite{CF85}, \cite{CF88}, \cite{TT88}})
\[
\frac{1}{2} \frac{d}{dt} \left|\chi_1(t)\wedge \cdots 
\wedge \chi_m(t) \right|_{\wedge^{m} L^2(\Om)}^2
= Tr \left(\widetilde{P}_m(t)\circ F^{\prime} (\widetilde{T}(t)) 
\circ \widetilde{P}_m(t)\right)  
\left|\chi_1(t)\wedge \cdots 
\wedge \chi_m(t) \right|_{\wedge^{m} L^2(\Om)}^2,
\]
which gives
\begin{equation}
\left|\chi_1(t)\wedge \cdots \wedge \chi_m(t)
\right|_{\wedge^{m} L^2(\Om)}^2 =
\left|\zeta_1 \wedge \cdots \wedge \zeta_m
\right|_{\wedge^{m} L^2(\Om)}^2 
\exp \int_0^t Tr \left(\widetilde{P}_m(s)\circ 
F^{\prime} (\widetilde{T}(s)) 
\circ \widetilde{P}_m(s)\right) ds,
\label{A-TRACE}
\end{equation}
where  $\chi_1(s), \cdots, \chi_m(s)$ are the solutions of 
(\ref{DIM-1})--(\ref{DW}) corresponding to the initial data
$\zeta_1, \cdots, \zeta_m,$ respectively. 
$Tr \left(\widetilde{P}_m(s)\circ F^{\prime} (\widetilde{T}(s))
\circ \widetilde{P}_m(s)\right) $ is the trace of the linear operator 
$\left(\widetilde{P}_m(s)\circ F^{\prime} (\widetilde{T}(s))
\circ \widetilde{P}_m(s) \right)$ and 
$\widetilde{P}_m(s)$ is the 
$L^2(\Om)$ orthogonal projector   
onto the space spanned by $\{ \chi_1(s), \cdots, \chi_m(s)\}.$
Thanks to (\ref{A-TRACE})
$\{\chi_1(s), \cdots, \chi_m(s)\}$ are linearly independent for every 
$s\geq 0$ if and only if $\{ \zeta_1, \cdots, \zeta_m\}$ are 
linearly independent. Hence, from now on we assume that 
$\{ \zeta_1, \cdots, \zeta_m\}$ are 
linearly independent.
Let $\{ \psi_1(s), \cdots, \psi_m(s) \}$ be an $L^2(\Om)$ 
orthonormal basis of the space 
spanned by $\{ \chi_1(s), \cdots, \chi_m(s)\}.$ Notice that 
$\{ \psi_1(s), \cdots, \psi_m(s) \}$ are in $V_2$ for $s>0.$
Thus we have  
\[
Tr \left(\widetilde{P}_m(s)\circ F^{\prime} (\widetilde{T}(s)) \circ
  \widetilde{P}_m(s)\right) 
= \sum_{j=1}^m \left( F^{\prime} (\widetilde{T}(s)) \psi_j(s) , 
\psi_j(s) \right).
\]
Observe that 
\begin{eqnarray*}
&&\hskip-.08in
\left( F^{\prime} (\widetilde{T}(s)) \psi_j(s) , \psi_j(s) \right)  
 = - \|\psi_j(s\|^2  + \int_{\Om} \left[  u_j
 \cdot \nabla (\widetilde{T}+T^*) 
+\widetilde{w}_j \pp_z \widetilde{T} \right] \psi_j(s) \; dxdydz,  
\end{eqnarray*}
where, for $j=1, 2, \cdots, m,$ 
$u_j (x,y,z,s)$ and $\widetilde{w}_j (x,y,z,s)$ 
are as in (\ref{DV-1})--(\ref{DW}), respectively,
but with $\psi_j(s)$ replacing $\chi$.
Following the same steps that led the estimates (\ref{W-4}),  we have
\begin{eqnarray*}
&&\hskip-.68in
\left|{ \int_{\Om}  u_j \cdot \nabla (\widetilde{T}+T^*) \;  
\psi_j(s)  }\right| \leq 
C \left( \| \widetilde{T}\|  + \|T^*\|_{H^1(M)} \right) 
 \; \left| \psi_j(s) \right| \; \left\| \psi_j(s) \right\|;
\end{eqnarray*}
and
\begin{eqnarray*} 
&&\hskip-.68in
\left|{ \widetilde{w}_j \pp_z \widetilde{T}\; \psi_j(s) }\right|
\leq  C \left| \widetilde{T} \right|^{\frac{3}{4}}
\left\| \nabla \widetilde{T} \right\|^{\frac{1}{4}}
\left| \psi_j(s) \right|^{\frac{1}{4}} 
\left\| \psi_j(s) \right\|^{\frac{7}{4}}.
\end{eqnarray*}
Recall that $|\psi_j|=1,$ for
$j=1, 2, \cdots, m.$ Thus, 
\begin{eqnarray*}
&&\hskip-.68in
\left| { \int_{\Om} \left[  u_j
 \cdot \nabla (\widetilde{T}+T^*) 
+\widetilde{w}  \pp_z \widetilde{T} \right] \psi_j(s) \; dxdydz } 
\right| \\
&&\hskip-.68in
\leq  C \left( \| \widetilde{T}\|  + \|T^*\|_{H^1(M)} \right) 
\; \left\| \psi_j(s) \right\| +   
C \left| \widetilde{T} \right|^{\frac{3}{4}}
\left\| \nabla \widetilde{T} \right\|^{\frac{1}{4}}
\left\| \psi_j(s) \right\|^{\frac{7}{4}}.
\end{eqnarray*} 
By using Young's inequality and the above estimate,  we have
\begin{eqnarray*}
&&\hskip-.68in
\left( F^{\prime} (\widetilde{T}(s)) \psi_j(s) , \psi_j(s) \right) 
\leq - \frac{1}{2} \|\psi_j(s\|^2 + 
C \left[\; \| \widetilde{T} (s)\|^2 + \|T^*\|^2_{H^1(M)} 
+ \left| \widetilde{T} \right|^6 
\left\| \nabla \widetilde{T} \right\|^2  \right].
\end{eqnarray*}
By (\ref{LL}), we have
\[
\sum_{j=1}^m \|\psi_j(s\|^2 \geq \ll_1 + \cdots + \ll_m 
\geq C \ll_1 m^2.
\]
As a result, we obtain
\[
Tr \left(\widetilde{P}_m(s)\circ F^{\prime} (\widetilde{T}(s)) 
\circ \widetilde{P}_m(s) \right)
\leq - C \ll_1 m^2 + 
C \left[\; \| \widetilde{T} (s)\|^2 + \|T^*\|^2_{H^1(M)} 
+ \left| \widetilde{T} \right|^6 
\left\| \nabla \widetilde{T} \right\|^2  \right].
\]
Hence,
\begin{eqnarray*}
&&\hskip-.68in
\frac{1}{t} \int_0^t Tr \left( \widetilde{P}_m(s)\circ 
F^{\prime} (\widetilde{T}(s)) 
\circ \widetilde{P}_m(s) \right) \; ds  \\
&&\hskip-.6in
\leq  
- C \ll_1 m^2 +C \left[\; \frac{1}{t} \int_0^t  \| \widetilde{T}
  (s)\|^2 \; ds + \|T^*\|^2_{H^1(M)} + 
\frac{1}{t} \int_0^t \left| \widetilde{T} \right|^6 
\left\| \nabla \widetilde{T} \right\|^2  \; ds \right].
\end{eqnarray*}
Therefore, by applying (\ref{L-W-I}), we get
\begin{eqnarray*}
&&\hskip-.68in
\limsup_{t\rightarrow\infty}\sup_{|\widetilde{T}_0| \leq R_a (T^*,Q)} 
\sup_{\begin{array}{c} 
     \zeta_j\in L^2(\Om)  \\
     |\zeta_j| \leq 1 \\
     j=1, \cdots, m
     \end{array} } 
\left[ \frac{1}{t} \int_0^t Tr \left( \widetilde{P}_m(s)\circ 
F^{\prime} (\widetilde{T}(s)) \right) 
\circ \widetilde{P}_m(s) \; ds 
\right] \\
&&\hskip-.68in
\leq -C \ll_1 m^2 + K_4 (T^*, Q),
\end{eqnarray*}
where
\begin{equation}
K_4 (T^*, Q)=  C (R_a(Q, T^*))^6 
\left[ 1 + \| T^*\|_{H^1 (M)}^2 + |Q|^2 \right].  \label{K-4}
\end{equation}
In order to guarantee 
$-C\ll_1 m^2 + K_4 (T^*, Q) 
\leq 0$ we need to
choose $m$ large enough such that
\begin{eqnarray*}
&&\hskip-.28in
m > C \left( \frac{K_4 (T^*, Q)}{\ll_1} \right)^{1/2}.
\end{eqnarray*}
Therefore, the Hausdorff and fractal dimensions of the attractor 
$\mathcal{A}$ can be estimated by (cf. e.g., {\bf \cite{EFNT}})
\[
d_H (\mathcal{A}) \leq d_F (\mathcal{A}) \leq  
C \left( \frac{K_4 (T^*, Q)}{\ll_1} \right)^{1/2}. 
\]
\end{proof}

\noindent
\section*{Acknowledgments}
This work was supported in part by the NSF
grants No. DMS--0204863 and DMS--0204794, by
the DOE under contract number W--7405--ENG--36, and by the US CRDF under grant number
RM1--2343--MO--02. This research was completed while E.S.T. was the Stanislaw M. Ulam 
Visiting Scholar at the CNLS in the Los Alamos National laboratory.


\begin{thebibliography}{99}


\bibitem{AR75}  R.A. Adams, {\em Sobolev Spaces},
Academic Press, New York, 1975.

\bibitem{CT01} C. Cao and E.S. Titi, {\em Global Well--posedness and 
Finite Dimensional Global Attractor for a 3--D Planetary Geostrophic 
Viscous Model},  Comm. Pure Appl. Math. {\bf 56} (2003), 198-233.

\bibitem{CJ55} J.G. Charney, {\em The gulf stream as an inertial
boundary layer,} Proc. Nat. Acad. Sci. U.S.A., {\bf 41}(1955), 731--740.

\bibitem{C84}  I. Chavel, {\em Eigenvalues in Riemannian Geometry}.
Academic Press, 1984.

\bibitem{CV86} A. Colin de Verdiere, {\em On mean flow
instability within the planetary geostrophic equations,}
J. Phys. Oceanogr., {\bf 16}, (1986), 1981--1984.



\bibitem{CF85} P. Constantin and C. Foias, {\em Global Lyapunov 
exponents, Kaplan--Yorke formulas and the dimension of the 
attractors for $2$D Navier--Stokes
equations}, Comm. Pure Appl. Math. {\bf 38} (1985), 1--27.


\bibitem{CF88} P. Constantin and C. Foias, {\em Navier-Stokes Equations,}
  The University of Chicago Press, 1988.


\bibitem{EFNT} A. Eden, C. Foias, B. Nicolaenko and R. Temam, {\em
Exponential Attractors for Dissipative Evolution Equations,}
Research in Applied Mathematics, {\bf 37}, Masson, Paris, 1994. 

\bibitem{FP84} P. Fabrie, {\em 
Solution forte tridimensionnelle d'un probl\'{e}me de 
convection--diffusion en milieu poreux. (French) [Strong
three-dimensional solutions for a convection--diffusion problem in a
porous medium]}, C. R. Acad. Sci. Paris S\'{e}r. I Math. 
{\bf 298} (1984), 249--251.



\bibitem{FP89} P. Fabrie, {\em Solutions fortes et majorations 
asymptotiques pour le mod\'{e}le de Darcy-Forchheimer en convection 
naturelle. (French)
[Strong solutions and asymptotic upper bounds for 
the Darcy-Forchheimer model in natural convection]},
 Ann. Fac. Sci. Toulouse Math. {\bf 10} (1989), 7--26. 


\bibitem{GA94} G.P. Galdi, {\em An Introduction to the Mathematical
Theory of the Navier-Stokes Equations,} Vol. I \& II,
Springer-Verlag, 1994.


\bibitem{Killworth}  P. D. Kilworth, {\em A two-level wind and buoyancy
driven thermocline model,} {\bf 15} (1985), 1414--1432. 
 
\bibitem{LADY} O.A. Ladyzhenskaya,  {\em The Boundary 
Value Problems of Mathematical Physics,} Springer-Verlag, 1985.

\bibitem{LADY91} O.A. Ladyzhenskaya, {\em Attractors for Semigroups 
and Evolution Equations,}  Cambridge University Press, Cambridge, 1991.


\bibitem{LL69} J.L. Lions, {\em Quelques M\'{e}thodes de R\'{e}solution 
des Probl\`{e}mes aux Limites
Non--lin\'{e}aires,}  Dunod, Paris, 1969.


\bibitem{LOT96} C.D. Levermore, M. Oliver and E.S. Titi, 
{\em Global well-posedness for models of shallow water in 
a basin with a varying bottom,} Indiana Univ. Math. J., 
{\bf 45}(1996), 479--510. 

\bibitem{PJ84} J. Pedlosky, {\em The equations for geostrophic
motion in the ocean,} J. Phys.  Oceanogr.,  
{\bf 14} (1984), 448--455.

\bibitem{PJ87} J. Pedlosky, {\em Geophysical Fluid Dynamics,}
Springer-Verlag, New York, 1987.

\bibitem{Ph} N. A. Phillips, Geostrophic motion, {\it Rev. Geophys.},
{\bf 1} (1963), 123--176. 

\bibitem{RS59} A.R. Robinson and H. Stommel, 
{\em The oceanic thermocline and the associated thermohaline
circulation,} Tellus, {\bf 11}(1959), 295--308.


\bibitem{Salmon-86} R. Salmon, {\em A simplified linear ocean
circulation theory,} J. Mar. Res., {\bf 44} (1986),
695--711.


\bibitem{Salmon-90} R. Salmon, {\em The thermocline as an ``internal
boundary layer'',} J. Mar. Res., {\bf 48} (1990), 437--469. 


\bibitem{SR97} R. Samelson, {\em Coastal boundary conditions 
and the baroclinic structure of wind--driven continental 
shelf currents,}
Journal of Physical Oceanography, {\bf 27}(1997), 2645--2662. 

\bibitem{SAMELSON} R. Samelson, private communication. 

\bibitem{STW98} R. Samelson, R. Temam and S. Wang,
{\em Some mathematical properties of the planetary geostrophic  equations
for large scale ocean circulation,} Applicable analysis, {\bf 70}(1998),
147--173.

\bibitem{STW00} R. Samelson, R. Temam and S. Wang,
{\em Remarks on the planetary geostrophic model of gyre scale ocean 
circulation}, Differential and Integral Equation, {\bf 13} (2000), 1--14.

\bibitem{SV97} R. Samelson and G.K. Vallis, 
{\em A simple friction and diffusion scheme for planetary geostrophic
basin models,}
Journal of Physical Oceanography, {\bf 27}(1997), 186--194. 

\bibitem{SV97A} R. Samelson and G.K. Vallis, 
{\em Large--scale circulation with small diapycnal diffusion: 
the two--thermocline limit,}
J. Mar. Res., {\bf 55}(1997), 223--275. 

\bibitem{SD96} D. Seidov, {\em An intermediate model for 
large--scale ocean circulation studies,}
Dynamics of Atmospheres and Oceans, {\bf 25}(1996), 25--55. 

\bibitem{SH48} H. Stommel, {\em The westward intensification of
wind--driven ocean currents,} Trans. Amer. Geophys. union, 
{\bf 29}(1948), 291--304.

\bibitem{TT79} R. Temam, {\em Navier-Stokes Equations, Theory and
Numerical Analysis,} North-Holland, 1979.


\bibitem{TT83}  R. Temam, {\em Navier-Stokes Equations and
Nonlinear Functional Analysis,} CBMS Regional Conference series,
No. {\bf 41}, SIAM, Philadelphia, 1983.

\bibitem{TT88} R. Temam, {\em Infinite-Dimensional Dynamical Systems in
Mechanics and Physics,} Applied Mathematical Sciences, {\bf 68},
Springer-Verlag, New York, 1988.


\bibitem{WE59} P. Welander, 
{\em An advective model of the  ocean thermocline,} Tellus, 
{\bf 11}(1959), 309--318.

\bibitem{WE71} P. Welander, 
{\em Some exact solutions to the equations describing an
ideal--fluid thermocline,} J. Mar. Res., 
{\bf 21}(1971), 60--68.


\bibitem{Winton-Sarachik} M.  Winton and E. Sarachik, {\em 
Thermocline
oscillations induced by strong steady salinity forcing of ocean
general circulation models,} J. Phys. Oceanogr.,
{\bf 23}  (1993), 1389--1410.

\bibitem{Zhang} S. Zhang, C. A. Lin, R. J. Greatbatch, {\em A thermocline
model for ocean-climate studies,} J. Mar. Res., {\bf 50} (1992),
99--124.

\end{thebibliography}
\end{document}